\documentclass[11pt,oneside,reqno]{amsart}

\usepackage[a4paper, total={450pt,675pt}]{geometry}
\usepackage[OT2,T1]{fontenc}
\usepackage[utf8]{inputenc}
\usepackage{amssymb,bm}

\usepackage[dvipsnames]{xcolor}
\usepackage[
    colorlinks=true,
    linkcolor=Maroon,
    citecolor=JungleGreen,
    urlcolor=NavyBlue]{hyperref}

\usepackage{xcolor,colortbl}
\usepackage[capitalize,nameinlink]{cleveref}

\usepackage{enumitem}
\usepackage{dsfont}
\usepackage{tikz}
\usetikzlibrary{patterns}
\usepackage{subcaption}
\usetikzlibrary{arrows}
\usepackage{multirow}

\usepackage{lipsum}
\usepackage[normalem]{ulem}
\newtheorem{theorem}{Theorem}[section]

\newtheorem{lemma}[theorem]{Lemma}
\newtheorem{proposition}[theorem]{Proposition}
\newtheorem{corollary}[theorem]{Corollary}

\newtheorem{remark}[theorem]{Remark}

\crefname{section}{Sect.}{section}
\numberwithin{equation}{section}
\DeclareMathOperator{\re}{Re}

\DeclareMathOperator{\im}{Im}
\DeclareMathOperator{\bigO}{O}

\DeclareMathOperator{\supp}{supp}

\DeclareMathOperator{\const}{const.}

\newcommand*\diff{\mathop{}\!\mathrm{d}}

\makeatletter
\def\@tocline#1#2#3#4#5#6#7{\relax
  \ifnum #1>\c@tocdepth % then omit
  \else
    \par \addpenalty\@secpenalty\addvspace{#2}%
    \begingroup \hyphenpenalty\@M
    \@ifempty{#4}{%
      \@tempdima\csname r@tocindent\number#1\endcsname\relax
    }{%
      \@tempdima#4\relax
    }%
    \parindent\z@ \leftskip#3\relax
    \advance\leftskip\@tempdima\relax
    \rightskip\@pnumwidth plus4em \parfillskip-\@pnumwidth
    #5\leavevmode\hskip-\@tempdima
      \ifcase #1
       \or\or \hskip 2em \or \hskip 2em \else \hskip 3em \fi%
      #6\nobreak\relax
    \dotfill\hbox to\@pnumwidth{\@tocpagenum{#7}}\par
    \nobreak
    \endgroup
  \fi}
\makeatother
\usepackage[
    backend=biber,
    style=alphabetic,
    giveninits=true,
    maxbibnames=10,
    sorting=nyt]{biblatex}

\DeclareLabelalphaTemplate{
  \labelelement{
    \field[final]{shorthand}
    \field{label}
    \field[strwidth=3]{labelname}}
  \labelelement{
    \field[strwidth=2,strside=right]{year}}
}

\DeclareFieldFormat[article,incollection]{title}{#1}
\DeclareFieldFormat{date}{#1}
\DeclareFieldFormat{pages}{{#1}}
\DeclareFieldFormat{doi}{
\href{https://www.doi.org/#1}{Doi}}
\DeclareFieldFormat{url}{
\href{#1}{Doi}}

\renewbibmacro{in:}{%
  \ifentrytype{article}{}{\printtext{\bibstring{in}\intitlepunct}}}

\begin{document}
%%%%%%%%%%%%%%%%%%%%%%%%%%%%%%%%%%%%%%%%%%%%%%%%%%%%%%%%%%%%
%%%%%%%%%%               TITLE PAGE               %%%%%%%%%%
%%%%%%%%%%%%%%%%%%%%%%%%%%%%%%%%%%%%%%%%%%%%%%%%%%%%%%%%%%%%
\title[Smoothing properties of dispersive equations
on non-compact symmetric spaces]{
Smoothing properties of dispersive equations \\
on non-compact symmetric spaces}

\author{Vishvesh Kumar, Michael Ruzhansky, and Hong-Wei Zhang}

\begin{abstract}
We establish the Kato-type smoothing property, i.e., global-in-time smoothing estimates with homogeneous weights, for the Schrödinger equation on Riemannian symmetric spaces of non-compact type and general rank. These form a rich class of manifolds with nonpositive sectional curvature and exponential volume growth at infinity, e.g., hyperbolic spaces. We achieve it by proving the Stein-Weiss inequality and the resolvent estimate of the corresponding Fourier multiplier, which are of independent interest. Moreover, we extend the comparison principles to symmetric spaces and deduce different types of smoothing properties for the wave equation, the Klein-Gordon equation, the relativistic and general orders Schrödinger equations. In particular, we observe that some smoothing properties, which are known to fail on the Euclidean plane, hold on the hyperbolic plane.

\end{abstract}

\keywords{dispersive equation, non-compact symmetric space, smoothing property, comparison principle, resolvent estimate, Stein-Weiss inequality}

%%%%%%%%%%% MSC2020
\makeatletter
\@namedef{subjclassname@2020}{\textnormal{2020}
    \it{Mathematics Subject Classification}}
\makeatother
\subjclass[2020]{22E30, 35J10, 35P25, 43A85, 43A90.}
%%%%%%%%%%% MSC2020

\maketitle
\tableofcontents

%%%%%%%%%%%%%%%%%%%%%%%%%%%%%%%%%%%%%%%%%%%%%%%%%%%%%%%%%%%%
%%%%%%%%%%               SECTION I                %%%%%%%%%%
%%%%%%%%%%%%%%%%%%%%%%%%%%%%%%%%%%%%%%%%%%%%%%%%%%%%%%%%%%%%
\section{Introduction}\label{Section.1 Introduction}
%%%%%%%%%%%%%%%%%%%%%%%%%%%%%%%%%%%%%%%%%%%%%%%%%%%%%%%%%%%%
A central topic in the area of analysis is to understand the influence of geometry on the behavior of solutions to the nonlinear evolution partial differential equations. Dispersive equations are well-known examples for their different phenomena on curved manifolds. In the study of nonlinear dispersive equations, there are two main tools: the \textit{Strichartz inequality} and the \textit{smoothing property}. Strichartz inequality usually means $L^{p}L^{q}$ space‐time mixed norm estimates, while smoothing properties are Sobolev $L^{2}$ space‐time estimates. Both approaches rely on global analysis for solutions to the corresponding linearized equations.

Non-compact symmetric spaces are Riemannian manifolds with nonpositive sectional curvature. Because of their exponential growth at infinity and the validity of the Kunze-Stein phenomenon, there are better dispersive properties and then stronger Strichartz inequalities. Such phenomenon was first observed in real hyperbolic spaces, which are the simplest models of non-compact symmetric spaces of rank one, see \cite{Fon97,Ion00a,Ban07,AnPi09,IoSt09,MeTa11,APV12,AnPi14}. The generalization for arbitrary ranks is recently achieved in \cite{AnZh23,AMPVZ23}. On the other hand, on general \textit{compact} manifolds, one has only local-in-time Strichartz inequalities with some loss of derivatives. See, for instance, \cite{Bou93a,Bou93b,BGT04,GePi10}.

Motivated by these interesting phenomena for the Strichartz inequality, we ask naturally that, as another primary tool in the study of dispersive equations, {\it how the smoothing property works on non-compact symmetric spaces?} The smoothing properties have been extensively studied in the Euclidean setting over the past three decades. See \cite{Sjo87,CoSa88,Veg88,KaYa89,BeDe91,Wat91,BeKl92,Sim92,KPV98,Sug98,Wal02,Hos03,Sug03,Chi08,RuSu12,DAn15,RuSu16} and references therein. We will discuss some of these works in detail in the following subsections by comparing them with our results on symmetric spaces. 

In the non-Euclidean backgrounds, there is also much literature on the \textit{local-in-time} smoothing properties. See, for instance, \cite{CKS95,Doi96,Doi00,Bur04,MMT08,Dat09,BGH10,ChWu13,ChMe14,BHS20}. The present paper focuses on the \textit{global-in-time} smoothing properties, which are less known than the rich local-in-time theory. As highlighted in \cite{RoTa07}, the main difficulty is that, besides the semiclassical analysis in high frequency, one also requires a more detailed analysis in low and medium frequencies. In that paper, Rodnianski and Tao established the global-in-time smoothing estimate with \textit{inhomogeneous weights} on asymptotically flat manifolds obeying the non-trapping condition. See also \cite{BoHa10,VaWu10,Bou11} for relevant estimates in this setting. Similar estimates with inhomogeneous weights were previously considered in symmetric spaces \cite{Kai14} and graded Lie groups \cite{Man17}. See also \cite{LLOS18,GeLe21} for recent progress on other related problems in hyperbolic spaces, such as smoothing estimates for the Schrödinger equation with potentials or $L^p$-estimates with $p>2$.

In this paper, we begin by establishing the Kato-type smoothing properties, namely, global-in-time smoothing estimates with \textit{homogeneous weights} for the Schrödinger equation. This is achieved by proving the resolvent estimate and the Stein-Weiss inequality, which are of independent interests as well. We emphasize that our setting does not enjoy the dilation property, hence the common rescaling method fails. Our second main result is to generalize the comparison principles from \cite{RuSu12} to non-compact symmetric spaces. This robust method allows us to deduce different types of smoothing properties for the wave equation, Klein-Gordon equation, and Schrödinger-type equations with general orders (even with some time-variable coefficients). In particular, we observe that some estimates which are known to fail on the Euclidean plane, hold on the hyperbolic plane. Most of our arguments rely on the harmonic analysis on Riemannian symmetric spaces.

For simplicity, from now on, a symmetric space always means the non-compact type, and the smoothing estimate refers to the global-in-time estimate. We will denote by $\Delta$ the Laplace-Beltrami operator on an $n$-dimensional symmetric space $\mathbb{X}$ and by $D$ its shifted Laplacian, see the following section for more details. To make the difference, let $\Delta_{\mathbb{R}^N}$ be the usual Laplacian in $\mathbb{R}^{N}$ ($N\ge2$) and $\widetilde{D}_{x}=(-i\partial{x_{1}},\,...\,,-i\partial{x_{N}})$. For $x\in\mathbb{X}$ or $x\in\mathbb{R}^{N}$, we denote by $|x|$ the (geodesic) distance between $x$ and the origin, and let $\langle{x}\rangle=(1+|x|^{2})^{1/2}$. Throughout the paper, the notation $a\lesssim{b}$ between two positive expressions means that $a\le{C}b$ for some constants $C>0$, and $a\asymp{b}$ means $a\lesssim{b}\lesssim{a}$.

\subsection{Smoothing estimates}
Consider the free Schrödinger equation in $\mathbb{R}^{N}$:
\begin{align}\label{S1 Sch} 
    (i\partial_{t}+\Delta_{\mathbb{R}^{N}})\,u(t,x)\,=\,0,
    \qquad\,u(0,x)\,=u_{0}(x),
\end{align}
whose solution is given by $u(t,x)=e^{it\Delta_{\mathbb{R}^{N}}}u_{0}(x)$. It is known that the solution operator $e^{it\Delta_{\mathbb{R}^{N}}}$ preserves the $L^2$-norm for each fixed time $t\in\mathbb{R}$. The smoothing property is a regularity improvement in the sense that, we could gain extra regularity (in comparison with the initial data) by integrating the solution to \eqref{S1 Sch} in time. More precisely, the solution to \eqref{S1 Sch} satisfies the smoothing property:
    \begin{align}\label{S1 BSmoothing}
        \|B(x,\widetilde{D}_{x})u\|_{
        L^{2}(\mathbb{R}_{t}\times\mathbb{R}^{N})}\,
        \lesssim\,\|u_{0}\|_{L^{2}(\mathbb{R}^{N})},
    \end{align}
where $B(x,\widetilde{D})$ is one of the following operators:

\begin{table}[ht]
\setlength{\tabcolsep}{20pt}
\renewcommand{\arraystretch}{2}
\begin{tabular}{|c|c|c|c|c|}
\hline
\cellcolor{gray!25} Type 
& \cellcolor{gray!25}  $B(x,\widetilde{D})$ 
& \cellcolor{gray!25}  Regularity condition \\
\hline
(I)
&$|x|^{\alpha-1}|\widetilde{D}|^{\alpha}$
& $1-\frac{N}{2}<\alpha<\frac12$ \\ 
\hline
(II)
& $\langle{x}\rangle^{-s}|\widetilde{D}|^{\frac12}$
& $s>\frac12$\\ 
\hline
(III)
&$\langle{x}\rangle^{-s}\langle{\widetilde{D}}\rangle^{\frac12}$
& $s\ge1$ ($s>1$ if $N=2$)  \\
\hline
\end{tabular}
\vspace{10pt}
\caption{Regularity conditions in $\mathbb{R}^{N}$ ($N\ge2$) for the Schrödinger equation.}
\label{S1 TableSchRN}
\end{table}

\textit{Restriction theorem} and \textit{resolvent estimate} (or their variants) are two main standard methods to deduce the above smoothing estimates. In \cite{RuSu12}, the authors introduced two other tools: the \textit{canonical transform} and the \textit{comparison principle}. The first helps to simplify the smoothing estimate to some 1-dimensional estimates, and the latter allows one to transfer smoothing estimates among different equations. In this paper, we will deduce the Kato-type smoothing property by establishing the resolvent estimate, and widen its regularity range by proving the Stein-Weiss inequality. By extending the comparison principle to symmetric spaces, we deduce new smoothing estimates for wave and Klein-Gordon equations.

%%%%%%%%%%%%%%%%%%%%%%%%%%%%%%%%%%%%%%%%%%%%%%%%%%%%%%%%%%%%%%%%%%%%%%%%%%%%%%%%%%
\subsection{Statement of main results}
Consider a non-compact symmetric space $\mathbb{X}=G/K$ of rank $\ell\ge1$, where $G$ and $K$ are suitable Lie groups. Let $n\ge2$ and $\nu\ge3$ be its manifold dimension and dimension at infinity (or pseudo-dimension, see Section \ref{Section.2 Prelim} for more details about these notations). Let $\Delta$ be the Laplace-Beltrami operator on $\mathbb{X}$ and denote by $D^{2}=-\Delta-|\rho|^{2}$ its shifted Laplacian. Here $|\rho|^{2}$ is the bottom of the $L^2$ spectrum of $-\Delta$ on $\mathbb{X}$. We consider the natural Schrödinger equation
    \begin{align}\label{S1 Schrodinger X}
    \begin{cases}
        i\partial_{t}u(t,x)\,+\,D_{x}^{2}u(t,x)\,=\,0,
        \qquad\,t\in\mathbb{R},\,\,\,x\in\mathbb{X},\\[5pt]
        u(0,x)\,=\,u_{0}(x),
    \end{cases}
    \end{align}
whose solution is given by $u(t,x)=e^{itD_{x}^{2}}u_{0}(x)$. The first part of this article focuses on the following smoothing property.

\begin{theorem}[Kato-type smoothing property]\label{main thm smoothing}
Let $\mathbb{X}$ be a symmetric space of dimension $n\ge3$ and pseudo-dimension $\nu\ge3$. Suppose that $1-\min\lbrace{\frac{n}{2},\frac{\nu}{2}}\rbrace<\alpha<\frac{1}{2}$. Then, the solution to the Schrödinger equation \eqref{S1 Schrodinger X} satisfies the smoothing property
    \begin{align}\label{main thm smoothing schrodinger}
        \||x|^{\alpha-1}\,D_{x}^{\alpha}\,u
        \|_{L^2(\mathbb{R}_{t}\times\mathbb{X})}\,
        \lesssim\,\|u_{0}\|_{L^2(\mathbb{X})}.
    \end{align}
Moreover, if $\mathbb{X}$ is of dimension $n=2$, then \eqref{main thm smoothing schrodinger} holds for all $-\frac{1}{2}<\alpha<\frac{1}{2}$.
\end{theorem}

\begin{remark}
The regularity condition in Theorem \ref{main thm smoothing} is optimal in some special cases. For example, when $G/K$ is a symmetric space with $G$ complex, the manifold dimension $n$ and the pseudo-dimension $\nu$ coincide, and the estimate \eqref{main thm smoothing schrodinger} fails for any $\alpha\le1-\frac{\nu}{2}$ or $\alpha\ge\frac12$, see Remark \ref{S3 remark optimal}.
\end{remark}

\begin{remark}\label{S1 Kato}
In Kato's theory, for a self-adjoint operator $H$ in a separable Hilbert space $\mathcal{H}$, one says that a densely-defined closed operator $A$ on $\mathcal{H}$ is $H$-smooth if 
   \begin{align*}
        |\im\big((H-\zeta)^{-1}\,A^{*}f,\,A^{*}f\big)|\,
        \lesssim\,
            \|f\|_{\mathcal{H}}^{2},
        \qquad\forall\,\zeta\in\mathbb{C}\smallsetminus{\mathbb{R}}.
    \end{align*}
Moreover, it is known that $A$ is $H$-smooth if and only if the smoothing property
    \begin{align*}
        \int_{\mathbb{R}}\diff{t}\,
        \|Ae^{-itH}f\|_{\mathcal{H}}^{2}\,
        \lesssim\,\|f\|_{\mathcal{H}}^{2}
    \end{align*}
holds, see \cite{Kat66,KaYa89}. Theorem \ref{main thm smoothing} is equivalent to say that, for suitable $\alpha$, the operator $|x|^{\alpha-1}D_{x}^{\alpha}$ is $D^{2}$-smooth on $\mathbb{X}$. In the $2$-dimensional case, $\mathbb{X}=\mathbb{H}^{2}$ is a hyperbolic plane. It has rank $\ell=1$ and pseudo-dimension $\nu=3$. Theorem \ref{main thm smoothing} shows that $|x|^{\alpha-1}D^{\alpha}$ is $D^{2}$-smooth on $\mathbb{H}^{2}$ for all $-\frac{1}{2}<\alpha<\frac{1}{2}$, while the operator  $|x|^{\alpha-1}\widetilde{D}^{\alpha}$ is $(-\Delta_{\mathbb{R}^{2}})$-smooth on $\mathbb{R}^{2}$ if and only if $0<\alpha<\frac{1}{2}$. It follows in particular that the weight $|x|^{-1}$ is $D^{2}$-smooth on $\mathbb{H}^{2}$.
\end{remark}

Note that the usual rescaling argument, which is used to establish the smoothing estimate with homogeneous weights in the Euclidean space, is not valid in the current setting. We require delicate and different analysis around or away from the origin. Our Theorem \ref{main thm smoothing} is achieved by combining the following resolvent estimate and the Stein-Weiss inequality.

\begin{theorem}[Resolvent estimate]\label{main thm resolv}
Let $\mathbb{X}$ be a symmetric space of rank $\ell\ge1$. Suppose that $-\frac{1}{2}<\alpha<\frac{1}{2}$ if $\ell=1$ and $1-\frac{\ell}{2}<\alpha<\frac{1}{2}$ if $\ell\ge2$. Then, for all $f\in{L^2}(\mathbb{X})$, we have
   \begin{align}\label{S3 resolvent}
        \sup_{\zeta\in\mathbb{C}\smallsetminus\mathbb{R}}\,
        |\im(D^{2\alpha}(D^{2}-\zeta)^{-1}\,f,\,f)_{L^{2}(\mathbb{X})}|\,
        \lesssim\,
            \||\cdot|^{1-\alpha}\,f\|_{L^{2}(\mathbb{X})}^{2}.
    \end{align}
\end{theorem}

On symmetric spaces of higher rank ($\ell\ge2$), Theorem \ref{main thm resolv} shows that the smoothing property \eqref{main thm smoothing schrodinger} holds when $1-\frac{\ell}{2}<\alpha<\frac12$. We widen this regularity range to $1-\min\lbrace{\frac{n}{2},\frac{\nu}{2}}\rbrace<\alpha<\frac{1}{2}$ by proving the Stein-Weiss inequality on symmetric spaces. This inequality is known as the Hardy-Littlewood-Sobolev inequality with double weights, see \cite{StWe58}. As an elementary tool in harmonic analysis, it has been extended to many other non-Euclidean settings, such as the Heisenberg group \cite{HLZ12}, the Carnot group \cite{GMS10}, and the homogeneous Lie group \cite{KRS19}. Recall that these three settings enjoy the dilation property.

On symmetric spaces, the $L^{p}$-$L^{q}$-boundedness of the operator $(-\Delta-|\rho|^{2}+\xi^{2})^{-\sigma/2}$, with $\xi\ge0$ and $\re\sigma\ge0$, were progressively established in the 90s, see for instance \cite{Str83,Var88,Loh89,Ank92,CGM93}. We refer to \cite[pp. 109-111]{CGM93} for a review of these works. See also \cite{BFL08,MaSa08,Bec15,LLY20,Bec21} for studies on the best constant problem in real hyperbolic spaces. The authors in \cite{KKR23} have recently established the Stein-Weiss inequality on symmetric spaces for the operator $(-\Delta-|\rho|^{2}+\xi^{2})^{-\sigma/2}$, with $\xi>0$ large enough. In that case, the corresponding convolution kernel provides additional exponential decay at infinity. Hence, their approach can not be applied in our limiting case where $\xi=0$, which needs a more delicate analysis. The following theorem is an $L^2$ generalization of the Stein-Weiss inequality associated with the operator $D^{-\sigma}=(-\Delta-|\rho|^{2})^{-\sigma/2}$ on symmetric spaces.

\begin{theorem}[Stein-Weiss inequality]\label{main thm SW}
Let $\mathbb{X}$ be a symmetric space of dimension $n\ge2$ and pseudo-dimension $\nu\ge3$. Suppose that $\sigma>0$ and $\gamma_{1},\gamma_{2}<\min\lbrace{\frac{n}{2},\frac{\nu}{2}}\rbrace$ satisfy $\sigma=\gamma_{1}+\gamma_{2}$. Then, for all $f\in{L^{2}}(\mathbb{X})$, we have
    \begin{align}\label{main thm SW ineq}
        \|\,|\cdot|^{-\gamma_{1}}\,D^{-\sigma}\,f\|_{L^{2}(\mathbb{X})}\,
        \lesssim\,
        \||\cdot|^{\gamma_{2}}\,f\|_{L^{2}(\mathbb{X})}.
    \end{align}
\end{theorem}

\begin{remark}
The condition $\gamma_{1}+\gamma_{2}=\sigma>0$ excludes the case where $\gamma_{1}=\gamma_{2}=0$. Recall that the operator $D^{-\sigma}$ is not $L^{2}$-bounded for any $\sigma>0$, see for instance \cite{Loh89}. When $\sigma=0$, the inequality \eqref{main thm SW ineq} holds if and only if $\gamma_{1}=\gamma_{2}=0$.
\end{remark}

\begin{remark}\label{S1 rmk SWoptimal}
Since the convolution kernel of $D^{-\sigma}$ behaves differently depending on whether it is close to or away from the origin, the manifold dimension and the pseudo-dimension both play essential roles, and the conditions $\gamma_{1},\gamma_{2}<\min\lbrace{\frac{n}{2},\frac{\nu}{2}}\rbrace$ in Theorem \ref{main thm SW} are both necessary, see Remark \ref{S3 rmk SWoptimal}. However, we observe from Theorem \ref{main thm resolv} that the smoothing property \eqref{main thm smoothing schrodinger} holds in rank one for all $1-\frac{\nu}{2}<\alpha<\frac{1}{2}$ ($\nu=3$). In other words, its regularity range depends only on the pseudo-dimension, which is related to the vanishing order of the Plancherel measure at the origin. We conjecture that it would be the same case in general ranks, but it is difficult to reach by relying only on the Stein-Weiss inequality.
\end{remark}

In the second part of this paper, we extend the comparison principles from \cite{RuSu12} to symmetric spaces, see Theorem \ref{S4 Comparison principle} and Corollary \ref{S4 Secondary comparison principle}. Based on the smoothing properties of the Schrödinger equation, we deduce different types of smoothing estimates for other equations. Let us start with the Schrödinger-type equations with general orders. Consider the Cauchy problem with order $m>0$:
    \begin{align}\label{S1 CPm}
        (i\partial_{t}+D_{x}^{m})\,u(t,x)=\,0,
        \qquad\,u(0,x)\,=\,u_{0}(x),
    \end{align}
whose solution is given by $u(t,x)=e^{itD_{x}^{m}}u_{0}(x)$. The following result describes all three types of smoothing properties on symmetric spaces for the Schrödinger-type equation \eqref{S1 CPm}.

\begin{theorem}\label{main thm SchSmoothing}
Let $\mathbb{X}$ be a symmetric space of rank $\ell\ge1$, dimension $n\ge2$ and pseudo-dimension $\nu\ge3$. Suppose that $m>0$ and $A(x,D)$ is defined as in Table \ref{S1 TableSch}. Then $A(x,D)$ is $D^{m}$-smooth on $\mathbb{X}$, namely, the solution to the Cauchy problem \eqref{S1 CPm} satisfies the smoothing property:
    \begin{align}\label{main SmoothSchA}
        \|A(x,D_{x})\,u\|_{
        L^{2}(\mathbb{R}_{t}\times\mathbb{X})}\,
        \lesssim\,
        \|u_{0}\|_{L^{2}(\mathbb{X})}.
    \end{align}
\end{theorem}

\begin{table}[ht]
\setlength{\tabcolsep}{5pt}
\renewcommand{\arraystretch}{2}
\begin{tabular}{|c|c|c|c|}
\hline
\cellcolor{gray!25} Type
& \cellcolor{gray!25} $A(x,D)$ 
& \cellcolor{gray!25} $\ell=1$ 
& \cellcolor{gray!25} $\ell\ge2$ \\
\hline
\textnormal{(I)}
& $|x|^{\alpha-\frac{m}{2}}D^{\alpha}$
& $\frac{m-3}{2}<\alpha<\frac{m-1}{2}$
& $\frac{m-\min\lbrace{n,\nu}\rbrace}{2}<\alpha<\frac{m-1}{2}$\\ 
\hline
\textnormal{(II)}
&$\langle{x}\rangle^{-s}D^{\frac{m-1}{2}}$
& \multicolumn{2}{c|}{$s>\frac{1}{2}$ \textnormal{and} $m>0$} \\  
\hline
\textnormal{(III)}
&$\langle{x}\rangle^{-s}\langle{D}\rangle^{\frac{m-1}{2}}$
& \multicolumn{2}{c|}{{$s\ge\frac{m}{2}$ \textnormal{and} $1<m<\nu$}} \\
\hline
\end{tabular}
\vspace{10pt}
\caption{Regularity conditions on $\mathbb{X}$ for the Schrödinger-type equations with order $m>0$.}
\label{S1 TableSch}
\vspace{-15pt}
\end{table}

\begin{remark}\label{S1 RemarkHn}
In particular, we can write the above smoothing estimates in hyperbolic spaces $\mathbb{H}^{n}$ as follows:
    \begin{align*}
        \||x|^{\alpha-\frac{m}{2}}\,D_{x}^{\alpha}\,u\|_{
        L^{2}(\mathbb{R}_{t}\times\mathbb{H}^{n})}\,
        &\lesssim\,
        \|u_{0}\|_{L^{2}(\mathbb{H}^{n})},
        \qquad\,\tfrac{m-3}{2}<\alpha<\tfrac{m-1}{2}\,\,\,\textnormal{and}\,\,\,m>0,
        \\[5pt]
        \|\langle{x}\rangle^{-s}\,D_{x}^{\frac{m-1}{2}}\,u\|_{
        L^{2}(\mathbb{R}_{t}\times\mathbb{H}^{n})}\,
        &\lesssim\,
        \|u_{0}\|_{L^{2}(\mathbb{H}^{n})},
        \qquad\,s>\tfrac{1}{2}\,\,\,\textnormal{and}\,\,\,m>0,
        \\[5pt]
        \|\langle{x}\rangle^{-s}\,\langle{D_{x}}\rangle^{\frac{m-1}{2}}\,u\|_{
        L^{2}(\mathbb{R}_{t}\times\mathbb{H}^{n})}\,
        &\lesssim\,
        \|u_{0}\|_{L^{2}(\mathbb{H}^{n})},
        \qquad\,s\ge\tfrac{m}{2}\,\,\,\textnormal{and}\,\,\,1<m<3,
    \end{align*}
where the regularity conditions depend only on the pseudo-dimension $\nu=3$.
\end{remark}

Let us compare the regularity conditions on symmetric spaces with the ones in the Euclidean setting. Recall that $B(x,\widetilde{D})$ is $\widetilde{D}_{x}^{m}$-smooth on $\mathbb{R}^{N}$ if it belongs to one of the following:

\begin{table}[ht]
\setlength{\tabcolsep}{20pt}
\renewcommand{\arraystretch}{2}
\begin{tabular}{|c|c|c|c|c|}
\hline
\cellcolor{gray!25} Type 
& \cellcolor{gray!25} $B(x,\widetilde{D})$ 
& \cellcolor{gray!25} Regularity condition in $\mathbb{R}^{N}$ \\
\hline
(I)
&$|x|^{\alpha-\frac{m}{2}}|\widetilde{D}|^{\alpha}$
& $\frac{m-N}{2}<\alpha<\frac{m-1}{2}\,\,\,\textnormal{and}\,\,\,m>0$ \\ 
\hline
(II)
& $\langle{x}\rangle^{-s}|\widetilde{D}|^{\frac{m-1}{2}}$
& $s>\frac{1}{2}$ and $m>0$  \\ 
\hline
(III)
&$\langle{x}\rangle^{-s}\langle{\widetilde{D}}\rangle^{\frac{m-1}{2}}$
& $s\ge\frac{m}{2}$ and $1<m<N$ ($s>\frac{m}{2}$ if $N=2$)\\
\hline
\end{tabular}
\vspace{10pt}
\caption{Regularity conditions on $\mathbb{R}^{N}$ ($N\ge2$) for the Schrödinger-type equations with order $m>0$.}
\vspace{-15pt}
\end{table}

For the Schrödinger equation ($m=2$) in $\mathbb{R}^{N}$, the Type (I) smoothing estimate was established by Kato and Yajima for $0\le{\alpha}<\frac12$ when $N\ge3$ and $0<{\alpha}<\frac12$ when $N=2$, see \cite{KaYa89}. Sugimoto \cite{Sug98} extended their regularity range to $1-\frac{N}{2}<\alpha<\frac12$ for all $N\ge2$, which is sharp and clarifies why $\alpha=0$ must be excluded when $N=2$. As a consequence, we know that $|x|^{-1}$ is not $(-\Delta_{\mathbb{R}^{2}})$-smooth on $\mathbb{R}^{2}$ as we mentioned in Remark \ref{S1 Kato}. In \cite[Theorem 5.2]{RuSu12}, the authors obtained the Type (I) estimate in $\mathbb{R}^{N}$ for all $N\ge2$ and $m>0$ satisfying $\frac{m-N}{2}<\alpha<\frac{m-1}{2}$, and pointed out that all the cases with different orders $m$ are equivalent to each other according to the comparison principle. This is also the case on symmetric spaces.

The Type (II) smoothing estimate in $\mathbb{R}^{N}$ ($N\ge2$) was proved in \cite{BeKl92} when $m=2$ and in \cite{Chi08} when $m>1$. A simpler proof based on the canonical transform was given in \cite[Theorem 5.1]{RuSu12}, where the authors showed that the type (II) estimate holds in fact for all $m>0$. By extending the arguments carried out in \cite{Chi08}, Kaizuka proved this estimate on symmetric spaces for $m>1$, see \cite{Kai14}. Using the comparison principle, we show that it holds for all $m>0$ as in the Euclidean setting. Note that $m=1$ is an important case corresponding to wave-type equations.

The Type (III) estimate has also been partially proved in \cite{Kai14}: the author showed that it holds on higher rank ($\ell\ge2$) symmetric spaces for all $1<m<\ell$. As a consequence of our improved inhomogeneous Stein-Weiss inequality (Corollary \ref{S3 SWCor}), the Type (III) estimate holds in fact on general symmetric spaces in the full range $1<m<\nu$. Note that this regularity condition is sharp and depends only on the pseudo-dimension $\nu$. In particular, it indicates that $\langle{x}\rangle^{-1}\langle{D}\rangle^{1/2}$ is $D^{2}$-smooth on $\mathbb{H}^{2}$. Recall that $\langle{x}\rangle^{-s}\langle{\widetilde{D}}\rangle^{1/2}$ is  $(-\Delta_{\mathbb{R}^2})$-smooth on $\mathbb{R}^{2}$ if and only if $s>1$, see \cite{KaYa89,Wal02, Chi08, RuSu12}.

Based on the above smoothing properties of Schrödinger-type equations and using comparison principles, we deduce smoothing estimates for some \textit{time-degenerate} and \textit{relativistic} Schrödinger equations as well, see Section \ref{subsection other examples}. Other noteworthy consequences of the comparison principles are the following smoothing properties of the wave and Klein-Gordon equations. 
    
\begin{theorem}\label{main smoothingKG}
Let $\mathbb{X}$ be a symmetric space of rank $\ell\ge1$, dimension $n\ge2$ and pseudo-dimension $\nu\ge3$. Consider the Cauchy problem
    \begin{align}\label{S1 KG}
        \begin{cases}
            (\partial_{t}^{2}+D_{x}^{2}+\zeta)\,u(t,x)=\,0,\\[5pt]
            u(0,x)\,=\,u_{0}(x),\,\partial_{t}|_{t=0}\,u(t,x)\,=\,u_{1}(x),
        \end{cases}
    \end{align}
which is the wave equation when $\zeta=0$ and the Klein-Gordon equation when $\zeta>0$. Let $s>\frac12$. Suppose that $-1<\beta<0$ when $\ell=1$ and $1-\min\lbrace{\frac{n}{2},\frac{\nu}{2}}\rbrace<\beta<\frac12$ when $\ell\ge2$. Then, the solution to the Cauchy problem \eqref{S1 KG} satisfies the following smoothing properties.
\begin{itemize}
    \item (Wave equation) If $\zeta=0$ , we have
    \begin{align}
        \|\langle{x}\rangle^{-s}\,u\|_{
        L^{2}(\mathbb{R}_{t}\times\mathbb{X})}\,
        &\lesssim\,
        \|u_{0}\|_{L^{2}(\mathbb{X})}\,+\,\|D_{x}^{-1}\,u_{1}\|_{L^{2}(\mathbb{X})},
        \label{S1 wave1}\\[5pt]
        \||x|^{\beta-\frac{1}{2}}\,D_{x}^{\beta}\,u\|_{
        L^{2}(\mathbb{R}_{t}\times\mathbb{X})}\,
        &\lesssim\,
        \|u_{0}\|_{L^{2}(\mathbb{X})}\,+\,\|D_{x}^{-1}\,u_{1}\|_{L^{2}(\mathbb{X})}.
        \label{S1 wave2}
    \end{align}
    
    \item (Klein-Gordon equation) If $\zeta>0$, we have
        \begin{align}\label{S1 wave3}
        \|\langle{x}\rangle^{-1}\,u\|_{
        L^{2}(\mathbb{R}_{t}\times\mathbb{X})}\,
        &\lesssim\,
        \|u_{0}\|_{L^{2}(\mathbb{X})}\,+\,\|D_{x}^{-1}\,u_{1}\|_{L^{2}(\mathbb{X})}.
    \end{align}
\end{itemize}
\end{theorem}

\begin{remark}\label{S1 RmkWave}
Smoothing properties of wave-type equations are well-known in the Euclidean setting, see for instance \cite{Ben94,RuSu12}. Since we deduce the above estimates from the Schrödinger equation, we observe different phenomena in $2$-dimensional cases as well. On the one hand, the estimate \eqref{S1 wave2} holds on $\mathbb{H}^{2}$ for all $-1<\beta<0$, while a similar estimate holds on $\mathbb{R}^{2}$ only for $-\frac12<\beta<0$. On the other hand, an estimate such as \eqref{S1 wave3} does not hold on $\mathbb{R}^2$ unless one considers the weight $\langle{x}\rangle^{-s}$ with $s>1$ instead of $\langle{x}\rangle^{-1}$.
\end{remark}

%%%%%%%%%%%%%%%%%%%%%%%%%%%%%%%%%%%%%%%%%%%%%%%%%%%%%%%%%%%%%%%%%%%%%%%%%%%%%%%% 
\subsection{Layout}
This paper is organized as follows. After a short review of harmonic analysis on symmetric spaces, we prove the Stein-Weiss inequality and establish the smoothing properties of the Schrödinger equation in Section \ref{Section.3 Smoothing}. In Section \ref{Section.4 Comparison}, we extend the comparison principles to symmetric spaces and deduce different types of smoothing properties for some other equations. Two technical lemmas are placed in Appendix.

%%%%%%%%%%%%%%%%%%%%%%%%%%%%%%%%%%%%%%%%%%%%%%%%%%%%%%%%%%%%%%%%%%%%%%%%%%%%%%%%%%
%%%%%%%%%%                            SECTION II                      %%%%%%%%%%%%
%%%%%%%%%%%%%%%%%%%%%%%%%%%%%%%%%%%%%%%%%%%%%%%%%%%%%%%%%%%%%%%%%%%%%%%%%%%%%%%%%%
\section{Preliminaries}\label{Section.2 Prelim}
In this section, we review briefly harmonic analysis on Riemannian symmetric spaces of non-compact type. We adopt the standard notation and refer to \cite{Hel78,Hel94,Hel00,GaVa88} for more details. 
    
\subsection{Non-compact symmetric spaces}
Let $G$ be a semisimple Lie group, connected, non-compact, with the finite center, and let $K$ be a maximal compact subgroup of $G$. The homogeneous space $\mathbb{X}=G/K$ is a Riemannian symmetric space of non-compact type. Let $\mathfrak{g}=\mathfrak{k}\oplus\mathfrak{p}$ be the Cartan decomposition of the Lie algebra $\mathfrak{g}$ of $G$. The Killing form of $\mathfrak{g}$ induces a $K$-invariant inner product $\langle.\,,\,.\rangle$ on $\mathfrak{p}$ and therefore a $G$-invariant Riemannian metric on $\mathbb{X}$.
%%%%%
Fix a maximal abelian subspace $\mathfrak{a}$ in $\mathfrak{p}$. The rank of $\mathbb{X}$ is the dimension $\ell$ of $\mathfrak{a}$. We identify $\mathfrak{a}$ with its dual $\mathfrak{a}^{*}$ by means of the inner product inherited from $\mathfrak{p}$. Let $\Sigma\subset\mathfrak{a}$ be the root system of $( \mathfrak{g},\mathfrak{a})$. Once a positive Weyl chamber $\mathfrak{a}^{+}\subset\mathfrak{a}$ has been selected, $\Sigma^{+}$ (resp. $\Sigma_{r}^{+}$ or $\Sigma_{s}^{+}$) denotes the corresponding set of positive roots (resp. positive reduced roots or simple roots).
%%%%%
Let $n$ be the dimension and $\nu$ be the dimension at infinity (or pseudo-dimension) of $\mathbb{X}$: 
    \begin{align}\label{S2 Dimensions}
        n\,=\,
        \ell\,+\,\sum_{\alpha \in \Sigma^{+}}\,m_{\alpha}
        \quad\textnormal{and}\quad
        \nu\,=\,\ell\,+\,2|\Sigma_{r}^{+}|,
    \end{align}
where $m_{\alpha}$ is the dimension of the positive root subspace $\mathfrak{g}_{\alpha}$. Notice that these two dimensions behave differently depending on the geometric structure of $\mathbb{X}$. For example, $\nu=3$ while $n\ge2$ is arbitrary in rank one, $\nu=n$ if $G$ is complex, and $\nu=2n-\ell>n$ when $G$ is split.

Let $\mathfrak{n}$ be the nilpotent Lie subalgebra of $\mathfrak{g}$ associated to $\Sigma^{+}$ and let $N=\exp\mathfrak{n}$ be the corresponding Lie subgroup of $G$. We have decompositions 
    \begin{align*}
        \begin{cases}
            \,G\,=\,N\,(\exp\mathfrak{a})\,K
            \quad&\textnormal{(Iwasawa)}, \\[5pt]
            \,G\,=\,K\,(\exp\overline{\mathfrak{a}^{+}})\,K
            \quad&\textnormal{(Cartan)}.
        \end{cases}
    \end{align*}
On the one hand, we can write the Haar measure $\diff{g}$ on $G$ in the Cartan decomposition:
    \begin{align*}
        \int_{G}\diff{g}\,f(g)\,
        =\,\const\,\int_{K}\,\diff{k}_{1}\,
        \int_{\mathfrak{a}^{+}}\,\diff{g}^{+}\,
        \underbrace{\prod_{\alpha\in\Sigma^{+}}\,
            (\sinh\langle{\alpha,g^{+}}\rangle)^{m_{\alpha}}
        }_{\delta(g^{+})}\, 
        \int_{K}\,\diff{k}_{2}\,f(k_{1}(\exp g^{+})k_{2}).
    \end{align*}
Notice that $\langle{\alpha,g^{+}}\rangle$ is nonnegative for every $\alpha\in\Sigma^{+}$ and all $g^{+}\in\overline{\mathfrak{a}^{+}}$. Let $\rho\in\mathfrak{a}^{+}$ be the half sum of all positive roots counted with their multiplicities:
    \begin{align*}
        \rho\,
        =\,\frac{1}{2}\,\sum_{\alpha\in\Sigma^{+}} 
            \,m_{\alpha}\,\alpha.
    \end{align*}
The density $\delta(g^{+})$ satisfies 
     \begin{align}\label{S2 density}
        \delta(g^{+})\,
        \asymp\,
            \prod_{\alpha\in\Sigma^{+}}
            \Big\lbrace 
            \frac{\langle\alpha,g^{+}\rangle}
            {1+\langle\alpha,g^{+}\rangle}
            \Big\rbrace^{m_{\alpha}}\,
            e^{\langle2\rho,g^{+}\rangle}\,
        \lesssim\,
            \begin{cases}
                |g^{+}|^{n-\ell}\,
                &\quad\textnormal{if}\;\;|g^{+}|\le\,1,\\[5pt]
                e^{\langle2\rho,g^{+}\rangle}\,
                &\quad\textnormal{for all}\;\;g^{+}\in\overline{\mathfrak{a}^{+}}.
            \end{cases}
    \end{align}
%%%%%
On the other hand, we can normalize the Haar measure $\diff{g}$ such as
    \begin{align}
        \int_{G}\diff{g}\,f(g)\,
        =\,
        \int_{N}\diff{n}\,\int_{\mathfrak{a}}\diff{A}\,
        e^{\langle{-2\rho,A}\rangle}\,
        \int_{K}\diff{k}\,f(n(\exp{A})k)
    \end{align}
where $A=A(g)$ is the unique $\mathfrak{a}$-component of $g$ in the Iwasawa decomposition.

\subsection{Harmonic analysis on symmetric spaces}
The harmonic analysis has been well developed on symmetric spaces. In the present paper, we shall need different types of transforms, such as the Helgason-Fourier transform $\mathcal{F}$, the Harish-Chandra transform $\mathcal{H}$, the Radon transform $\mathcal{R}$, and the modified Radon transform $\mathcal{JR}$. We review their definitions and basic properties in the following. Bearing in mind that the Cartan subspace $\mathfrak{a}$ is an $\ell$-dimensional flat submanifold of $\mathbb{X}$, we denote by $\mathcal{F}_{\mathfrak{a}}$ the classical Fourier transform:
    \begin{align*}
        \mathcal{F}_{\mathfrak{a}}f(\lambda)\,
        =\,\int_{\mathfrak{a}}\diff{H}\,e^{-i\langle{\lambda,H}\rangle}f(H)\,
        \qquad\textnormal{and}\qquad
         \mathcal{F}_{\mathfrak{a}}^{-1}g(H)\,
        =\,\int_{\mathfrak{a}}\diff{\lambda}\,e^{i\langle{\lambda,H}\rangle}g(\lambda)\,
    \end{align*}
for suitable functions $f$ and $g$ on $\mathfrak{a}$.

%%%%%%%%%%%%%%%%%%%%%%%%%%%%%%%%%%%%%%%%%%%%%%%%%%%%%%%%%%%%%%%%%%%%%%%%%%%%%%%%%%
\subsubsection{Helgason-Fourier transform}
Let $f$ be a Schwarz function on $\mathbb{X}$. Denote by $M$ the centralizer of $\exp\mathfrak{a}$ in $K$ and $\diff{b}$ a $K$-invariant normalized measure on $B=K/M$. The Helgason-Fourier transform and its inverse formula are defined by
    \begin{align}\label{S2 Helgason Fourier}
        \mathcal{F}f(\lambda,kM)\,
        =\,
        \int_{G}\diff{g}\,
        e^{\langle{-i\lambda+\rho,A(k^{-1}g)}\,\rangle}\,f(g)
        \qquad\forall\,\lambda\in\mathfrak{a},\;\;
        \forall\,k\in{K},
    \end{align}
and
    \begin{align}\label{S2 Inverse Helgason Fourier}
        f(gK)\,=\,|W|^{-1}\,
        \int_{\mathfrak{a}}\diff{\lambda}\,|\mathbf{c}(\lambda)|^{-2}\,
        \int_{B}\diff{b}\,e^{\langle{i\lambda+\rho,\,A(k^{-1}g)}\rangle}\,
        \mathcal{F}f(\lambda,b),
        \qquad\forall\,g\in{G},
    \end{align}
where $W$ denotes the Weyl group associated to $\Sigma$, see for instance \cite[Ch.III]{Hel94}.
%%%%%
Here $|\mathbf{c}(\lambda)|^{-2}$ is the so-called \textit{Plancherel density} and can be expressed via the Gindikin-Karpelevič formula:
    \begin{align}\label{S2 Plancherel Density}
        |\mathbf{c}(\lambda)|^{-2}\,
        =\,\prod_{\alpha\in\Sigma_{r}^{+}}\,
            |\mathbf{c}_{\alpha}
            (\langle\alpha,\lambda\rangle)|^{-2},
    \end{align}
where each \textit{Plancherel factor} $|\mathbf{c}_{\alpha}(\cdot)|^{-2}$ is explicitly defined by some Gamma functions and extends to an analytic function in a neighbourhood of the real axis, see \cite[Theorem 6.14, p.447]{Hel00}. Moreover, for every $\alpha\in\Sigma_{r}^{+}$, the factor $|\mathbf{c}_{\alpha}(\cdot)|^{-2}$ is a homogeneous differential symbol of order $m_{\alpha}+m_{2\alpha}$ and satisfies
    \begin{align}
        |\mathbf{c}_{\alpha}(r)|^{-2}\,\asymp\,
        |r|^{2}(1+|r|)^{m_{\alpha}+m_{2\alpha}-2}\,
        \qquad\forall\,r\in\mathbb{R}.
    \end{align}
As a product of one-dimensional symbols, the Plancherel density $|\mathbf{c}(\lambda)|^{-2}$ is not a symbol on $\mathfrak{a}$ in general:
    \begin{align}\label{S2 Plancherel Estimates}
        \begin{cases}
        |\mathbf{c}(\lambda)|^{-2}\,
        \lesssim\,
        |\lambda|^{\nu-\ell} 
        & \textnormal{if \,} |\lambda|\le1,\\[5pt]
        |\nabla_{\mathfrak{a}}^{k}|\mathbf{c}(\lambda)|^{-2}|\,
        \lesssim\,
        |\lambda|^{n-\ell} 
        & \textnormal{if \,} |\lambda|\ge1
        \,\,\,\textnormal{and}\,\,\,k\in\mathbb{N}.
        \end{cases}
    \end{align}
The Plancherel theorem states that the Helgason-Fourier transform $\mathcal{F}$ extends to an isometry of $L^{2}(\mathbb{X})$ into $L^{2}(\mathfrak{a}\times{B},|W|^{-1}|\mathbf{c}(\lambda)|^{-2}\diff{\lambda}\diff{b})$, see for instance \cite[Theorem 1.5, p.227]{Hel94}.

%%%%%%%%%%%%%%%%%%%%%%%%%%%%%%%%%%%%%%%%%%%%%%%%%%%%%%%%%%%%%%%%%%%%%%%%%%%%%%%%%%
\subsubsection{Harish-Chandra transform}
A function $f$ is called bi-$K$-invariant on $\mathbb{X}$ if $f(k_{1}gk_{2})=f(g)$ for all $k_{1},k_{2}\in{K}$ and $g\in{G}$. With such functions, the Helgason-Fourier transform \eqref{S2 Helgason Fourier} reduces to
the Harish-Chandra transform
    \begin{align}
        \mathcal{H}f(\lambda)\,
        =\,\int_{G}\diff{g}\,\varphi_{-\lambda}(g)\,f(g) 
        \quad\forall\,\lambda\in\mathfrak{a},\;\;
        \forall\,f\in\mathcal{S}(K\backslash{G/K}),
        \label{S22 HarishChandra}
    \end{align}
where
    \begin{align*}
        \varphi_{\lambda}(g) 
        = \int_{K}\diff{k}\,
            e^{\langle{i\lambda+\rho,\,A(kg)}\rangle}
    \end{align*}
is the \textit{elementary spherical function}, see \cite[Theorem 4.3, p.418]{Hel00}. For every $\lambda$ in $\mathfrak{a}$, the spherical function $\varphi_{\lambda}$ is bi-$K$-invariant and satisfies $|\varphi_{\lambda}|\le\varphi_{0}$, where
    \begin{align}\label{S2 phi0}
        \varphi_{0}(\exp{H})\,
        \asymp\,\Big\lbrace{
        \prod_{\alpha\in\Sigma_{r}^{+}}(1+\langle{\alpha,H\rangle})\,
        }\Big\rbrace\,e^{-\langle{\rho,H}\rangle}
        \qquad\forall\,H\in\overline{\mathfrak{a}^{+}}.
    \end{align}

Denote by $\mathcal{S}(\mathfrak{a})^{W}$ the subspace of $W$-invariant functions in the Schwartz space $\mathcal{S}(\mathfrak{a})$. Then $\mathcal{H}$ is an isomorphism between $\mathcal{S}(K\backslash{G/K})$ and $\mathcal{S}(\mathfrak{a})^{W}$. The inverse formula of the Harish-Chandra transform is given by
    \begin{align}\label{Inverse Harish-Chandra}
        f(g)\,
        =\,\const\,\int_{\mathfrak{a}}\diff{\lambda}\,
            |\mathbf{c(\lambda)}|^{-2}\,
            \varphi_{\lambda}(g)\,
            \mathcal{H}f(\lambda) 
        \quad \forall\,g\in{G},\ 
            \forall\,f\in\mathcal{S}(K\backslash{G/K}).
    \end{align}

%%%%%%%%%%%%%%%%%%%%%%%%%%%%%%%%%%%%%%%%%%%%%%%%%%%%%%%%%%%%%%%%%%%%%%%%%%%%%%%%%%
\subsubsection{Radon and modified Radon transforms}
The Helgason-Fourier transform $\mathcal{F}$ is the flat Fourier transform $\mathcal{F}_{\mathfrak{a}}$ of the Radon transform
    \begin{align*}
        \mathcal{R}f(H,kM)\,=\,
        e^{-\langle{\rho,H}\rangle}\,
        \int_{N}\diff{n}\,f(n(\exp{H})k)
        \qquad\forall\,(H,kM)\in\,\mathfrak{a}\times{B}.
    \end{align*}
In other words, for all $\lambda\in\mathfrak{a}$ and $k\in{K}$, we can write
    \begin{align}\label{S2 transformXtoR}
        \mathcal{F}f(\lambda,kM)\,=\,
        \mathcal{F}_{\mathfrak{a}}[ \mathcal{R}f(\cdot,kM)](\lambda)
    \end{align}
see \cite[Ch.II, \S3]{Hel94}. This transform allows us to borrow some techniques from the Euclidean Fourier analysis in studying harmonic analysis on symmetric spaces. Denote by $\mathcal{J}$ the Fourier multiplier on $\mathfrak{a}$ with the symbol $|\mathbf{c}(\lambda)|^{-1}$. The $\mathcal{JR}$-transform is an isometry from $L^{2}(\mathbb{X})$ to $L^{2}(\mathfrak{a}\times{B},\,|W|^{-1}\diff{H}\diff{b})$.
Indeed, by using the Plancherel formula with respect to $\mathcal{F}_{\mathfrak{a}}$ and $\mathcal{F}$, we have
    \begin{align*}
        \|\mathcal{JR}f
        \|_{L^{2}(\mathfrak{a}\times{B},\,|W|^{-1}\diff{H}\diff{b})}\,
        &=\,
            \||\mathbf{c}(\lambda)|^{-1}\,\mathcal{F}_{\mathfrak{a}}\mathcal{R}f
            \|_{L^{2}(\mathfrak{a}\times{B},\,
                |W|^{-1}\diff{\lambda}\diff{b})}\\[5pt]
        &=\,
            \|\mathcal{F}f
            \|_{L^{2}(\mathfrak{a}\times{B},\,
                |W|^{-1}\,|\mathbf{c}(\lambda)|^{-2}\diff{\lambda}\diff{b})}\\[5pt]
        &=\,\|f\|_{L^{2}(\mathbb{X})}.
    \end{align*}
In \cite[Proposition 3.1]{Kai14}, the author proved that, for any $\sigma>0$, the $\mathcal{JR}$ transform is a continuous map from $L^{2}(\mathbb{X},\langle{x}\rangle^{2\sigma}\diff{x})$ to $L^{2}(\mathfrak{a}\times{B},\langle{H}\rangle^{2\sigma}\diff{H}\diff{b})$. The following lemma shows that similar continuity remains valid with homogeneous weights. Its proof is not so different from the original one, we include the details in Appendix for the sake of completeness.

\begin{lemma}\label{S2 JR continuity}
For any $\sigma\ge0$, we have
    \begin{align}
        \|\mathcal{JR}f\|_{L^{2}(\mathfrak{a}\times{B},\,
            |H|^{2\sigma}\diff{H}\diff{b})}\,
        \lesssim\,
        \|f\|_{L^{2}(\mathbb{X},\,|x|^{2\sigma}\diff{x})}.
    \end{align}
\end{lemma}

%%%%%%%%%%%%%%%%%%%%%%%%%%%%%%%%%%%%%%%%%%%%%%%%%%%%%%%%%%%%%%%%%%%%%%%%%%%%%%%%%%
%%%%%%%%%%                           SECTION III                      %%%%%%%%%%%%
%%%%%%%%%%%%%%%%%%%%%%%%%%%%%%%%%%%%%%%%%%%%%%%%%%%%%%%%%%%%%%%%%%%%%%%%%%%%%%%%%%
\section{Kato-type smoothing property on symmetric spaces}\label{Section.3 Smoothing}
In this section, we establish the Kato-type smoothing property for the Schrödinger equation on $\mathbb{X}$. We start with the Stein-Weiss inequality, namely, Theorem \ref{main thm SW}. Together with the resolvent estimate Theorem \ref{main thm resolv}, we deduce Theorem \ref{main thm smoothing}.

%%%%%%%%%%%%%%%%%%%%%%%%%%%%%%%%%%%%%%%%%%%%%%%%%%%%%%%%%%%%%%%%%%%%%%%%%%%%%%%%%%
\subsection{Stein-Weiss inequality}
Recall the operator $D^{-\sigma}=(-\Delta-|\rho|^2)^{-\sigma/2}$ where $\sigma>0$. Notice that the balance condition $\sigma=\gamma_{1}+\gamma_{2}$ with $\gamma_{1},\gamma_{2}<\min\lbrace{\frac{n}{2},\frac{\nu}{2}}\rbrace$ implies  $\sigma<\min\lbrace{n,\nu}\rbrace$. Denote by $k_{\sigma}$ the bi-$K$-invariant convolution kernel of the operator $D^{-\sigma}$, it satisfies 
    \begin{align}\label{S3 Riesz ker estim}
        k_{\sigma}(x)\,\asymp\,
        \begin{cases}
            |x|^{\sigma-n}
            &\qquad\textrm{if}\,\,\,|x|\le1\,\,\,\textnormal{and}\,\,\,
                0<\sigma<n,
            \\[5pt]
            |x|^{\sigma-\nu}\,\varphi_{0}(x)\,
            &\qquad\textrm{if}\,\,\,|x|\ge1\,\,\,\textnormal{and}\,\,\,
                0<\sigma<\nu,
        \end{cases}
    \end{align}
see \cite[Theorem 4.2.2]{AnJi99}. Then Theorem \ref{main thm SW} is equivalent to the following proposition.

\begin{proposition}\label{S3 SW prop}
Let $\mathbb{X}$ be a symmetric space of dimension $n\ge2$ and pseudo-dimension $\nu\ge3$. Let $\sigma>0$, $\gamma_{1},\gamma_{2}<\min\lbrace{\frac{n}{2},\frac{\nu}{2}}\rbrace$, and $\sigma=\gamma_{1}+\gamma_{2}$. Then the operator $\mathcal{T}$ defined by
    \begin{align}
        \mathcal{T}f(x)\,=\,
        \int_{\mathbb{X}}\diff{y}\,
            |x|^{-\gamma_{1}}\,k_{\sigma}(y^{-1}x)\,
            |y|^{-\gamma_{2}}\,f(y)
    \label{S32 opT}
    \end{align}
is bounded from $L^{2}(\mathbb{X})$ into $L^{2}(\mathbb{X})$.
\end{proposition}

Because of the contrasting behaviors of the convolution kernel, as well as the density volume, we shall need different arguments depending on whether $|x|$ and $|y|$ are small or large. Roughly saying, when $|x|$ and $|y|$ are both small, the density volumes grow polynomially. In this case, we extend a key lemma from \cite{StWe58} to symmetric spaces. If $|x|$ or $|y|$ is large, the density volume grows exponentially, and we need to combine the dyadic decomposition with suitable Kunze-Stein phenomena. The proof of Proposition \ref{S3 SW prop} is based on the following two lemmas.

\begin{lemma}[Kunze-Stein phenomena]\label{S3 KS lemma}
Let $g$ be a bi-$K$-invariant function in $\mathcal{S}(\mathbb{X})$. Then, for any $f\in\mathcal{S}(\mathbb{X})$, we have
    \begin{align}
        \|g*f\|_{L^2(\mathbb{X})}\,
        \lesssim\,
        \|g\varphi_{0}\|_{L^1(\mathbb{X})}\,\|f\|_{L^2(\mathbb{X})},
        \label{S3 KS1}
    \end{align}   
and
    \begin{align}
        \|g*f\|_{L^2(\mathbb{X})}\,
        \lesssim\,
        \|g\|_{L^2(\mathbb{X},\,\langle{x}\rangle^{\nu}\diff{x})}\,
        \|f\|_{L^2(\mathbb{X})}.
        \label{S3 KS2}
    \end{align}
\end{lemma}

\begin{remark}\label{rmk KS}
These properties are consequences of Herz's principle (see \cite{Her70}) with bi-$K$-invariant functions. The inequality \eqref{S3 KS1} follows by \cite{Cow97}, see also \cite{APV11,Zha21}. The weighted version \eqref{S3 KS2} was previously stated in \cite[Corollary 2.6]{Kai14} for $g\in{L^2(\mathbb{X},\,\langle{x}\rangle^{\sigma}\diff{x})}$ with $\sigma>\nu$. Our lemma shows that such inequality remains valid in the critical case where $\sigma=\nu$, provided that $g$ is bi-$K$-invariant. In fact, the function $g$ involves only the bi-$K$-invariant convolution kernel in our proof of the Stein-Weiss inequality, and this endpoint improvement is crucial.
\end{remark}

\begin{proof}
According to Remark \ref{rmk KS} and property \eqref{S3 KS1}, it is sufficient to show that
    \begin{align*}
        \|g\varphi_{0}\|_{L^1(\mathbb{X})}\,
        \lesssim\,
        \|g\|_{L^2(\mathbb{X},\,\langle{x}\rangle^{\nu}\diff{x})}.
    \end{align*}
Let $\chi_{0}\in\mathcal{C}_{c}^{\infty}(\mathbb{R}_{+})$ be a cut-off function such that $\supp\chi_{0}\subset[0,1]$ and $\chi_{0}=1$ on $[0,\frac12]$. For all $x\in\mathbb{X}$, we denote by $\widetilde{\chi}_{0}(x)=\chi_{0}(|x|)$ and 
    \begin{align*}
        \widetilde{\chi}_{j}(x)\,
        =\,\chi_{0}(2^{-j}|x|)\,-\,\chi_{0}(2^{-j+1}|x|)
        \qquad\forall\,j\ge1,
    \end{align*}    
which are all bi-$K$-invariant cut-off functions on $\mathbb{X}$. For every $j\ge1$, $\widetilde{\chi}_{j}$ is compactly supported in $\lbrace{x\in\mathbb{X}\,|\,2^{j-2}\le|x|\le2^{j}}\rbrace$. In particular, we have $\sum_{j\in\mathbb{N}}\widetilde{\chi}_{j}=1$ and 
    \begin{align*}
        \|g\|_{L^2(\mathbb{X},\,\langle{x}\rangle^{\nu}\diff{x})}\,
        \asymp\,
            \Big\lbrace{\sum_{j\in\mathbb{N}}\,2^{\nu{j}}\,
            \|\widetilde{\chi}_{j}^{1/2}g\|_{L^{2}(\mathbb{X})}^{2}
            }\Big\rbrace^{1/2}
        \qquad\forall\,g\in{\mathcal{S}(\mathbb{X})}.
    \end{align*}
By using the partition of unity and the Cauchy-Schwarz inequality, we have
\begin{align*}
    \|g\varphi_{0}\|_{L^1(\mathbb{X})}^{2}\,
    &\lesssim\,
    \sum_{j\in\mathbb{N}}\,
    \Big\lbrace{
    \int_{\mathbb{X}}\diff{x}\,\widetilde{\chi}_{j}(x)\,|g(x)|\,\varphi_{0}(x)
    }\Big\rbrace^{2}\\[5pt]
    &\le\,
    \sum_{j\in\mathbb{N}}\,
    \Big\lbrace{
    \int_{\mathbb{X}}\diff{x}\,\widetilde{\chi}_{j}(x)\,|g(x)|^{2}
    }\Big\rbrace\,
    \Big\lbrace{
    \int_{\mathbb{X}}\diff{x}\,\widetilde{\chi}_{j}(x)\,\varphi_{0}^{2}(x)
    }\Big\rbrace.
\end{align*}
According to the density estimate \eqref{S2 density} and estimate \eqref{S2 phi0} of the ground spherical function, we obtain
\begin{align*}
    \int_{\mathbb{X}}\diff{x}\,
    (\widetilde{\chi}_{0}(x)+\widetilde{\chi}_{1}(x))\,
    \varphi_{0}^{2}(x)\,
    <\,+\infty
\end{align*}
and
\begin{align*}
    \int_{\mathbb{X}}\diff{x}\,\widetilde{\chi}_{j}(x)\,\varphi_{0}^{2}(x)\,
    &\lesssim\,
    \int_{1\le|x^{+}|\le2^{j}}\diff{x^{+}}\,\delta(x^{+})\,
    \varphi_{0}^{2}(\exp{x^{+}})\\[5pt]
    &\lesssim\,
    \int_{1}^{2^{j}}\diff{r}\,r^{\nu-\ell+\ell-1}\,
    \lesssim\,2^{\nu{j}}
\end{align*}
for all $j\ge2$. We deduce that
    \begin{align*}
    \|g\varphi_{0}\|_{L^1(\mathbb{X})}^{2}\,
    \lesssim\,
    \sum_{j\in\mathbb{N}}\,
    2^{\nu{j}}\,\|\widetilde{\chi}_{j}^{1/2}g\|_{L^{2}(\mathbb{X})}^{2}\,
    \asymp\,
    \|g\|_{L^2(\mathbb{X},\,\langle{x}\rangle^{\nu}\diff{x})}^{2},
    \end{align*}
which completes the proof of Lemma \ref{S3 KS lemma}.
\end{proof}

The next lemma extends \cite[Lemma 2.1]{StWe58} to symmetric spaces. Notice that the additional condition \eqref{S3 SWlemma cdt} appears naturally, since the growth of volume on symmetric spaces has different behaviors. We will include its detailed proof in Appendix.

\begin{lemma}[Stein-Weiss lemma]\label{S3 SW lemma}
Let $\mathcal{K}:\mathbb{R}_{+}\times\mathbb{R}_{+}\rightarrow\mathbb{R}_{+}$ be a homogeneous function of degree $-\ell$ such that 
    \begin{align*}
        \int_{0}^{+\infty}\diff{s}\,s^{\tfrac{\ell}{2}-1}\,\mathcal{K}(1,s)\,<\,+\infty.
    \end{align*}
    Let $\kappa_{1}$, $\kappa_{2}$ be two bi-$K$-invariant functions on $G$ such that
    \begin{align}\label{S3 SWlemma cdt}
        \kappa_{j}(\exp{x}^{+})\,\delta^{1/2}(x^{+})\,
        =\,\bigO(1)
        \qquad\forall\,x\in{G},\,\,
        \forall\,j=1,\,2,
    \end{align}
    where $x^{+}\in\mathfrak{a}$ is the middle component of $x$ in the Cartan decomposition. Then the operator $S:L^{2}(G)\rightarrow{L^{2}}(G)$ defined by
    \begin{align*}
        Sf(x)\,
        =\,\kappa_{1}(x)\,
            \int_{G}\diff{y}\,
            \mathcal{K}(|x|,|y|)\,
            \kappa_{2}(y)\,f(y)
    \end{align*}
    is bounded.
\end{lemma}

Now, let us turn to the proof of Proposition \ref{S3 SW prop}.

\begin{proof}[Proof of Proposition \ref{S3 SW prop}]
We prove the $L^2$-boundedness of $\mathcal{T}$ in different cases depending whether $|x|$, $|y|$ and $|y^{-1}x|$ are small or large. For $k=1,2,...,7$, we define
    \begin{align}\label{S3 operators Tj}
    \mathcal{T}_{k}f(x)\,=\,
        \int_{\mathbb{X}}\diff{y}\,\psi_{k}(x,y)\,
            |x|^{-\gamma_{1}}\,k_{\sigma}(y^{-1}x)\,|y|^{-\gamma_{2}}\,f(y),
    \end{align}
where $\psi_{j}(x,y)$ are suitable cut-off functions which will be specified in each case. Recall that $\chi_{0}$ is a cut-off function compactly supported in $[0,1]$, and $\chi_{0}=1$ on $[0,\frac12]$. Denote by $\chi_{\infty}=1-\chi_{0}$.

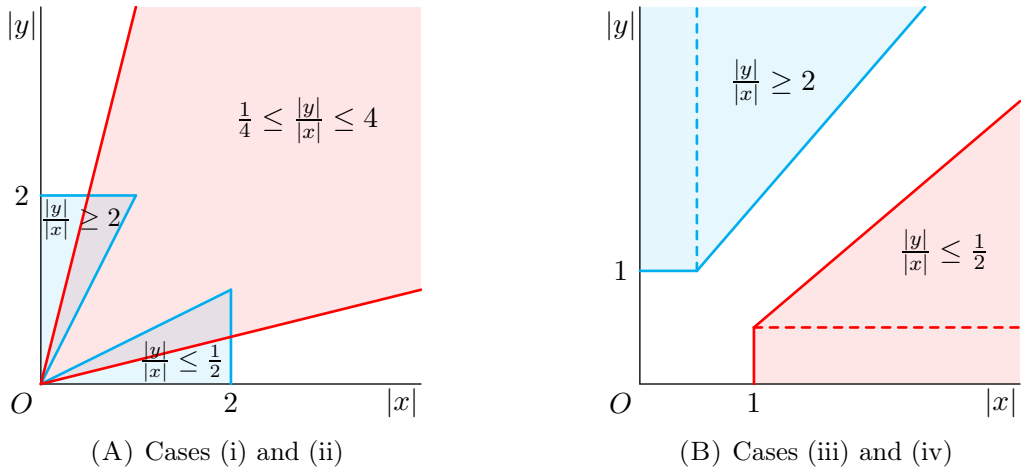
\begin{figure}[b]
     \begin{subfigure}[b]{0.49\textwidth}
     \centering
        \begin{tikzpicture}[line cap=round,line join=round,>=triangle 45,x=1.0cm,y=1.0cm, scale=0.5]
\clip(-1,-1) rectangle (10,10);
\fill[line width=0.pt,color=red,fill=red,fill opacity=0.10000000149011612] (2.5,10.) -- (10.,10.) -- (10.,2.5) -- (0.,0.) -- cycle;
\fill[line width=0.pt,color=cyan,fill=cyan,fill opacity=0.10000000149011612] (0.,0.) -- (0.,5.) -- (2.5,5) --cycle;
\fill[line width=0.pt,color=cyan,fill=cyan,fill opacity=0.10000000149011612] (0.,0.) -- (5,0) -- (5,2.5) --cycle;
\draw [line width=0.5pt] (0,0)-- (0,10);
\draw [line width=0.5pt] (0.,0.)-- (10,0);
%\draw [line width=1.pt,dashed] (10,7.5)-- (0,0);
%\draw [line width=1.pt,dashed] (0,0)-- (7.5,10.);
\draw [line width=1.pt,color=cyan] (0.,5)-- (2.5,5);
\draw [line width=1.pt,color=cyan] (5.,0.)-- (5,2.5);
\draw [line width=1.pt,color=cyan] (0,0)-- (2.5,5);
\draw [line width=1.pt,color=cyan] (0,0)-- (5,2.5);
\draw [line width=1.pt,color=red] (10.,2.5)-- (0.,0.);
\draw [line width=1.pt,color=red] (0,0)-- (2.5,10.);
\draw (-.5,-.5) node{$O$};
\draw (5,-.5) node{$2$};
\draw (9.5,-.5) node{$|x|$};
\draw (-.5,5) node{$2$};
\draw (-.5,9.5) node{$|y|$};
\draw (7,7) node{$\tfrac{1}{4}\le\tfrac{|y|}{|x|}\le4$};
\draw (1.05,4.35) node{\small $\tfrac{|y|}{|x|}\ge2$};
\draw (3.7,0.6) node{\small $\tfrac{|y|}{|x|}\le\frac12$};
\end{tikzpicture}
        \caption{Cases (i) and (ii)}
        \label{Fig Cases12}
     \end{subfigure}
     %\hfill
     \begin{subfigure}[b]{0.49\textwidth}
     \centering
        \begin{tikzpicture}[line cap=round,line join=round,>=triangle 45,x=1.0cm,y=1.0cm, scale=0.5]
\clip(-1,-1) rectangle (10,10);
\fill[line width=0.pt,color=red,fill=red,fill opacity=0.10000000149011612] (3.,0.) -- (10.,0.) -- (10.,7.5) -- (3.,1.5) -- cycle;
\fill[line width=0.pt,color=cyan,fill=cyan,fill opacity=0.10000000149011612] (0.,10.) -- (0.,3.) -- (1.5,3.) -- (7.5,10) -- cycle;
\fill[line width=0.pt,color=cyan,fill=cyan,fill opacity=0.10000000149011612,pattern=north west lines] (0,3) rectangle (1.5,10);
\fill[line width=0.pt,color=red,fill=red,fill opacity=0.10000000149011612,pattern=north west lines] (3,0) rectangle (10,1.5);
\draw [line width=0.5pt] (0,0)-- (0,10);
\draw [line width=0.5pt] (0.,0.)-- (10,0);
%\draw [line width=1.pt,dashed] (10.,7.5)-- (0,0);
%\draw [line width=1.pt,dashed] (0,0)-- (7.5,10.);
\draw [line width=1.pt,color=red] (3.,1.5)-- (3.,0.);
\draw [line width=1.pt,color=cyan] (0.,3.)-- (1.5,3.);
\draw [line width=1.pt,color=cyan,dashed] (1.5,3.)-- (1.5,10.);
\draw [line width=1.pt,color=red,dashed] (3.,1.5)-- (10,1.5);
\draw [line width=1.pt,color=red] (10.,7.5)-- (3,1.5);
\draw [line width=1.pt,color=cyan] (1.5,3)-- (7.5,10.);
\draw (-.5,-.5) node{$O$};
\draw (3,-.5) node{$1$};
\draw (9.5,-.5) node{$|x|$};
\draw (-.5,3) node{$1$};
\draw (-.5,9.5) node{$|y|$};
\draw (3.5,8) node{$\tfrac{|y|}{|x|}\ge2$};
\draw (8,3.5) node{$\tfrac{|y|}{|x|}\le\tfrac12$};
\end{tikzpicture}
        \caption{Cases (iii) and (iv)}
        \label{Fig Case3}
     \end{subfigure}
     \caption{Different cases in the study of $\mathcal{T}$.}
     \vspace{-20pt}
\end{figure}
    
\noindent\textbf{Case (i): if $|x|$ and $|y|$ are comparable.}
Let $\psi_{1}(x,y)=\chi_{0}(\frac{|y|}{4|x|})\chi_{\infty}(\frac{2|y|}{|x|})$. Then $\psi_{1}$ does not vanish when $\frac{1}{4}\le\frac{|y|}{|x|}\le4$ (see the red zone in Figure \ref{Fig Cases12}). We decompose dyadically
    \begin{align}\label{S3 dyadic}
        \|\mathcal{T}_{1}f\|_{L^{2}(\mathbb{X})}^{2}\,
        =\,
        \sum_{j\in\mathbb{Z}}\,
        \int_{\lbrace{x\in\mathbb{X}\,:\,
            2^{j}\le{|x|}\le2^{j+1}}\rbrace}\diff{x}\,
        \Big|{
        \int_{\mathbb{X}}\diff{y}\,\psi_{1}(x,y)\,
            |x|^{-\gamma_{1}}\,k_{\sigma}(y^{-1}x)\,
            |y|^{-\gamma_{2}}\,f(y)
        }\Big|^{2}.
    \end{align}    
Notice that $\frac{1}{4}\le\frac{|y|}{|x|}\le4$ and $2^{j}\le{|x|}\le2^{j+1}$ imply $2^{j-2}\le{|y|}\le2^{j+3}$ and $|y^{-1}x|\le2^{j+4}$. Then we have $|x|^{-\gamma_{1}}|y|^{-\gamma_{2}}\lesssim2^{-\sigma{j}}$, provided that $\gamma_{1}+\gamma_{2}=\sigma>0$. Hence
    \begin{align}\label{S3 case1 estim1}
        \|\mathcal{T}_{1}f\|_{L^{2}(\mathbb{X})}^{2}\,
        &\lesssim\,
        \sum_{j\in\mathbb{Z}}\,2^{-2\sigma{j}}\,
        \int_{\mathbb{X}}\diff{x}\,
        \Big|{
            \int_{\mathbb{X}}\diff{y}\,
            \chi_{0}(\tfrac{|y|}{2^{j+4}})\,
            \chi_{\infty}(\tfrac{|y|}{2^{j-2}})\,
            \chi_{0}(\tfrac{|y^{-1}x|}{2^{j+5}})\,k_{\sigma}(y^{-1}x)\,f(y)
        }\Big|^{2}\notag\\[5pt]
        &=\,\sum_{j\in\mathbb{Z}}\,2^{-2\sigma{j}}\,
            \big\|\big(\chi_{0}(\tfrac{|\cdot|}{2^{j+4}})\,
            \chi_{\infty}(\tfrac{|\cdot|}{2^{j-2}})\,
            f\big)
            *\big(\chi_{0}(\tfrac{|\cdot|}{2^{j+5}})k_{\sigma}\big)
            \big\|_{L^{2}(\mathbb{X})}^{2}.
    \end{align}
we deduce from the Kunze-Stein phenomenon \eqref{S3 KS1} that 
    \begin{align}\label{S3 case1 estim2}
        \big\|\big(\chi_{0}(\tfrac{|\cdot|}{2^{j+4}})\,
            \chi_{\infty}(\tfrac{|\cdot|}{2^{j-2}})\,
            f\big)
            &*\big(\chi_{0}(\tfrac{|\cdot|}{2^{j+5}})k_{\sigma}\big)
        \big\|_{L^{2}(\mathbb{X})}\notag\\[5pt]
        &\lesssim\,
        \big\|\chi_{0}(\tfrac{|\cdot|}{2^{j+4}})\,
            \chi_{\infty}(\tfrac{|\cdot|}{2^{j-2}})\,
            f\big\|_{L^{2}(\mathbb{X})}\,
        \big\|\chi_{0}(\tfrac{|\cdot|}{2^{j+5}})\,k_{\sigma}\,
            \varphi_{0}\big\|_{L^{1}(\mathbb{X})}.
    \end{align}
    
According to the kernel estimate \eqref{S3 Riesz ker estim}, the density estimate \eqref{S2 density}, and estimate \eqref{S2 phi0} of the ground spherical function, we obtain, on the one hand, 
    \begin{align}\label{S3 case1 estim3}
        \|\chi_{0}(\tfrac{|\cdot|}{2^{j+5}})\,k_{\sigma}\,
            \varphi_{0}\|_{L^{1}(\mathbb{X})}\,
        &\lesssim\,
            \int_{|x^{+}|\le1}\,\diff{x}^{+}\,
            |x^{+}|^{\sigma-n}\,|x^{+}|^{n-\ell}\,
            +\,
            \int_{1\le|x^{+}|\le2^{j+5}}\,\diff{x}^{+}\,
            |x^{+}|^{\sigma-\nu}\,
            |x^{+}|^{\nu-\ell}\notag\\[5pt]
        &=\,
            \int_{0}^{2^{j+5}}\,\diff{r}\,r^{\sigma-1}\,
            \lesssim\,2^{\sigma{j}}
    \end{align}
provided that $\sigma>0$ and $j\ge {-4}$.  On the other hand, if $j<{-4}$, we have similarly
    \begin{align}\label{S3 case1 estim4}
        \|\chi_{0}(\tfrac{|\cdot|}{2^{j+5}})\,k_{\sigma}\,
            \varphi_{0}\|_{L^{1}(\mathbb{X})}\,
        &\lesssim\,
            \int_{|x^{+}|\le2^{j+5}}\,\diff{x^{+}}\,
            |x^{+}|^{\sigma-n}\,|x^{+}|^{n-\ell}\notag\\[5pt]
        &=\,\int_{0}^{2^{j+5}}\,\diff{r}\,r^{\sigma-1}\,
        \lesssim\,2^{\sigma{j}}
    \end{align}
by using \eqref{S3 Riesz ker estim} and \eqref{S2 density} again.

We deduce from \eqref{S3 case1 estim1}, \eqref{S3 case1 estim2}, \eqref{S3 case1 estim3}, and \eqref{S3 case1 estim4}, that 
    \begin{align}\label{S3 estim T1}
        \|\mathcal{T}_{1}f\|_{L^{2}(\mathbb{X})}^{2}\,
        \lesssim\,
        \sum_{j\in\mathbb{Z}}\,
        \int_{\lbrace{x\in\mathbb{X}\,:\,
            2^{j}\le{|x|}\le2^{j+1}}\rbrace}\diff{x}\,
            |f(x)|^{2}\,
        =\,\|f\|_{L^{2}(\mathbb{X})}^{2}
    \end{align}
provided that $\gamma_{1}+\gamma_{2}=\sigma>0$.

\noindent\textbf{Case (ii): if $|x|$ and $|y|$ are not comparable, but both small.}
We denote by $\psi_{2}(x,y)=\chi_{0}(\frac{|x|}{2})\chi_{0}(\frac{2|y|}{|x|})$ and $\psi_{3}(x,y)=\chi_{0}(\frac{|y|}{2})\chi_{\infty}(\frac{|y|}{4|x|})$ two cut-off functions. Notice that $\supp\psi_{2}\cup\supp\psi_{3}$ corresponds to the blue zones in Figure \ref{Fig Cases12}. Since $|y^{-1}x|$ is also small in this case, we know from \eqref{S3 Riesz ker estim} that,
    \begin{align*}
        \mathcal{T}_{k}f(x)\,
        \asymp\,
        \int_{\mathbb{X}}\diff{y}\,
            \psi_{k}(x,y)\,|x|^{-\gamma_{1}}\,
            |y^{-1}x|^{\sigma-n}\,|y|^{-\gamma_{2}}\,f(y)
        \qquad\textnormal{with}\,\,\,k=2,\,3.
    \end{align*}

On the one hand, when $\frac{|y|}{|x|}\le\frac{1}{2}$, we have $|y^{-1}x|\ge\frac{|x|}{2}$ and then
    \begin{align*}
        \mathcal{T}_{2}f(x)\,
        &\lesssim\,
        \underbrace{\vphantom{\int}
        \chi_{0}(\tfrac{|x|}{2})\,
        |x|^{\tfrac{\ell-n}{2}}
        }_{\kappa_{1}(x)}\,
        \int_{\mathbb{X}}\diff{y}\,
        \underbrace{\vphantom{\int}
        \chi_{0}(\tfrac{2|y|}{|x|})\,
        |x|^{-\gamma_{1}+\sigma-n+\tfrac{n-\ell}{2}}\,
        |y|^{-\gamma_{2}+\tfrac{n-\ell}{2}}
        }_{\mathcal{K}(|x|,|y|)}\,
        \underbrace{\vphantom{\int}
        \chi_{0}(\tfrac{|y|}{2})\,
        |y|^{\tfrac{\ell-n}{2}}
        }_{\kappa_{2}(y)}
        f(y)\, 
    \end{align*}
provided that $0<\sigma<n$. Here, we observe that $(-\gamma_{1}+\sigma-n+\frac{n-\ell}{2})+(-\gamma_{2}+\frac{n-\ell}{2})=-\ell$, since $\gamma_1+\gamma_2=\sigma$, and
    \begin{align*}
        \chi_{0}(\tfrac{|z^{+}|}{2})\,|z^{+}|^{\tfrac{\ell-n}{2}}\,
        \delta(z^{+})^{\tfrac{1}{2}}\,
        =\,\bigO(1)
        \qquad\forall\,z\in\mathbb{X},
    \end{align*}
according to \eqref{S2 density}. Moreover, we have
    \begin{align*}
        \int_{0}^{\frac{1}{2}}\diff{s}\,
            s^{\tfrac{\ell}{2}-1}\,
            s^{-\gamma_{2}+\tfrac{n-\ell}{2}}\,
            <\,+\infty
    \end{align*}
provided that $\gamma_2<\frac{n}{2}$. We deduce from Lemma \ref{S3 SW lemma} that $\mathcal{T}_{2}$ is $L^{2}$-bounded. On the other hand, when $\frac{|y|}{|x|}\ge2$, we have $|y^{-1}x|\ge\frac{|y|}{2}$ and similarly
    \begin{align*}
        \mathcal{T}_{3}f(x)\,
        \lesssim\,
        \chi_{0}(\tfrac{|x|}{2})\,
        |x|^{\tfrac{\ell-n}{2}}\,
        \int_{\mathbb{X}}\diff{y}\,
        \chi_{\infty}({|y|}/{4|x|})\,
        |x|^{-\gamma_{1}+\tfrac{n-\ell}{2}}\,
        |y|^{-\gamma_{2}+\sigma-n+\tfrac{n-\ell}{2}}\,
        \chi_{0}(\tfrac{|y|}{2})\,
        |y|^{\tfrac{\ell-n}{2}}\,f(y).
    \end{align*}
We deduce from the same argument that $\mathcal{T}_{3}$ is $L^2$-bounded since
    \begin{align*}
        \int_{2}^{\infty}\diff{s}\,
            s^{\tfrac{\ell}{2}-1}\,
            s^{\gamma_{1}-n+\tfrac{n-\ell}{2}}\,
            <\,\infty
    \end{align*}
provided that $\gamma_{1}<\frac{n}{2}$.
    
Therefore, for all $\gamma_{1},\gamma_{2}<\frac{n}{2}$ such that $\gamma_{1}+\gamma_{2}= \sigma>0$, we have
    \begin{align}\label{S3 estim T2T3}
        \|\mathcal{T}_{k}f\|_{L^{2}(\mathbb{X})}\,
        \lesssim\,\|f\|_{L^{2}(\mathbb{X})}
        \qquad\textnormal{with}\,\,\,k=2,\,3.
    \end{align}

\noindent\textbf{Case (iii): if $|x|$ and $|y|$ are not comparable, and not both small or both large.}
We define $\psi_{4}(x,y)=\chi_{\infty}(\frac{|y|}{2})\chi_{0}(2|x|)$ and $\psi_{5}(x,y)=\chi_{\infty}(\frac{|x|}{2})\chi_{0}(2|y|)$. The support of $\psi_{4}$ (resp. $\psi_{5}$) corresponds to the blue (resp. red) shaded rectangle in Figure \ref{Fig Case3}. According to the local Harnack inequality of the elementary spherical function, we have $\varphi_{0}(y^{-1}x)\lesssim\varphi_{0}(y)$ when $|x|$ is bounded, see for instance \cite[Proposition 4.6.3]{GaVa88} or \cite[Remark 4.5]{APZ23}. Hence, for all $(x,y)\in\supp\psi_{4}$ and $0<\sigma<\nu$,
    \begin{align}\label{S3 phi0comp}
        k_{\sigma}(y^{-1}x)\,
        &\asymp\,
        |y^{-1}x|^{\sigma-\nu}\,
        \varphi_{0}(y^{-1}x)\notag\\[5pt]
        &\lesssim\,
        |y|^{\sigma-\nu}\varphi_{0}(y)\,
        \lesssim\,
        |y^{+}|^{\frac{2\sigma-\nu-\ell}{2}}\,e^{-\langle{\rho,y^{+}}\rangle}.
    \end{align}
    
Combining this estimate with the density estimates \eqref{S2 density}, we obtain
    \begin{align*}
        \|\mathcal{T}_{4}f(x)\|_{L^{2}(\mathbb{X})}^{2}\,
        &=\,
            \int_{\mathbb{X}}\diff{x}\,
            \Big|
            \int_{\mathbb{X}}\diff{y}\,
            \psi_{4}(x,y)\,|x|^{-\gamma_{1}}\,k_{\sigma}(y^{-1}x)\,
            |y|^{-\gamma_{2}}\,f(y)
            \Big|^{2}\\[5pt]
        &\lesssim\,
            \int_{\mathbb{X}}\diff{x}\,|x|^{-2\gamma_{1}}\,
            \Big\lbrace{
            \int_{\mathbb{X}}\diff{y}\,
            \psi_{4}^{2}(x,y)\,k_{\sigma}^{2}(y^{-1}x)\,|y|^{-2\gamma_{2}}
            }\Big\rbrace\,
            \|f\|_{L^{2}(\mathbb{X})}^{2},
    \end{align*}
where
    \begin{align*}
        \int_{\mathbb{X}}\diff{y}\,
            \psi_{4}^{2}(x,y)\,k_{\sigma}^{2}(y^{-1}x)\,|y|^{-2\gamma_{2}}\,
        &\lesssim\,\chi_{0}(2|x|)\,
        \int_{\mathfrak{a}^{+}}\diff{y^{+}}\,\delta(y^{+})\,
            \chi_{\infty}(\tfrac{|y^{+}|}{2})\,
            |y^{+}|^{2\sigma-\nu-\ell-2\gamma_{2}}\,e^{-2\langle{\rho,y^{+}}\rangle}\\[5pt]
        &\lesssim\,\chi_{0}(2|x|)
            \underbrace{
            \int_{1}^{+\infty}\diff{r}\,r^{2\gamma_{1}-\nu-1}
            }_{<\,+\infty},
    \end{align*}
provided that $\gamma_{1}+\gamma_{2}\ge\sigma$ and $\gamma_{1}<\frac{\nu}{2}$. On the other hand, we have
    \begin{align*}
        \int_{\mathbb{X}}\diff{x}\,\chi_{0}(2|x|)\,|x|^{-2\gamma_{1}}\,
        \lesssim\,
        \int_{|x^{+}|\le\frac12}\diff{x^{+}}\,
        |x^{+}|^{-2\gamma_{1}}\,|x^{+}|^{n-\ell}\,<\,+\infty,
    \end{align*}
provided that $\gamma_{1}<\frac{n}{2}$. Therefore, $\mathcal{T}_{4}$ is $L^{2}$-bounded. One can show the $L^2$-boundedness of $\mathcal{T}_{5}$ by using similar arguments. 

Therefore, for all $\gamma_{1},\gamma_{2}<\min\lbrace{\frac{n}{2},\frac{\nu}{2}}\rbrace$ satisfying $\gamma_{1}+\gamma_{2}\ge\sigma>0$, we have
    \begin{align}\label{S3 estim T4T5}
        \|\mathcal{T}_{k}f\|_{L^{2}(\mathbb{X})}\,
        \lesssim\,\|f\|_{L^{2}(\mathbb{X})}
        \qquad\textnormal{with}\,\,\,k=4,\,5.
    \end{align}

\noindent\textbf{Case (iv): if $|x|$ and $|y|$ are not comparable, but both large.}
In the last case, let us define $\psi_{6}(x,y)=\chi_{\infty}(\frac{|y|}{4|x|})\,\chi_{\infty}(|x|)$ and $\psi_{7}(x,y)=\chi_{0}(\frac{2|y|}{|x|})\,\chi_{\infty}(|y|)$, see the non-shaded blue and red triangles in Figure \ref{Fig Case3} for zones corresponding to their supports. For any $(x,y)\in\supp\psi_{6}$, we have $\frac{|y|}{|x|}\ge2$ and $|x|\ge\frac12$, which imply that $|y|\ge1$, $|x|\asymp\langle{x}\rangle$, and $\frac{|y|}{2}<|y^{-1}x|<\frac{3|y|}{2}$. Then, for any $\gamma_{1},\gamma_{2}<\frac{\nu}{2}$ satisfying $\gamma_{1}+\gamma_{2}\ge\sigma>0$, and for all $(x,y)$ in the support of $\psi_{6}$, we have
    \begin{align*}
        |x|^{-\gamma_{1}}\,|y|^{-\gamma_{2}}&= \langle{x}\rangle^{-\frac{\nu}{2}}|y^{-1}x|^{\frac{\nu}{2}-\sigma} 
        \underbrace{\vphantom{\int}
        |x|^{-\gamma_1} \langle x \rangle^{\frac{\nu}{2}}|y^{-1}x|^{\sigma-\frac{\nu}{2}} |y|^{-\gamma_2}
        }_{\lesssim\,|y|^{-\gamma_1+\frac{\nu}{2}+\sigma-\frac{\nu}{2}-\gamma_2}\,\le\,1}.
    \end{align*}
Combining this inequality with the partition of unity defined in the proof of Lemma \ref{S3 KS lemma}, we obtain
    \begin{align*}
        \|\mathcal{T}_{6}f(x)\|_{L^{2}(\mathbb{X})}^{2}\,
        &=\,
            \int_{\mathbb{X}}\diff{x}\,
            \Big|
            \int_{\mathbb{X}}\diff{y}\,\psi_{6}(x,y)\,
            |x|^{-\gamma_{1}}\,k_{\sigma}(y^{-1}x)\,
            |y|^{-\gamma_{2}}\,f(y)
            \Big|^{2}\\[5pt]
        &\lesssim\,
            \int_{\mathbb{X}}\diff{x}\,\langle{x}\rangle^{-\nu}\,
            \Big|
            \int_{\mathbb{X}}\diff{y}\,
            \psi_{6}(x,y)\,k_{\sigma}(y^{-1}x)\,
            |y^{-1}x|^{\frac{\nu}{2}-\sigma}\,f(y)
            \Big|^{2}\\[5pt]
        &\lesssim\,
            \sum_{j\in\mathbb{N}}\,
            \int_{\mathbb{X}}\diff{x}\,
            \langle{x}\rangle^{-\nu}\,
            \Big|
            \int_{\mathbb{X}}\diff{y}\,\widetilde{\chi}_{j}(y)\,
            \psi_{6}(x,y)\,k_{\sigma}(y^{-1}x)\,
            |y^{-1}x|^{\frac{\nu}{2}-\sigma}\,f(y)
            \Big|^{2}.
    \end{align*}

Recall that for all $y\in\supp\widetilde{\chi}_{j}\subset[2^{j-2},2^{j}]$ and $(x,y)\in\supp\psi_{6}$, we have $\chi_{0}(\frac{|y^{-1}x|}{2^{j+2}})=1$ and $\chi_{\infty}(2|y^{-1}x|)=1$. We deduce from the duality of \eqref{S3 KS2} that
    \begin{align*}
        \|\mathcal{T}_{6}f(x)\|_{L^{2}(\mathbb{X})}^{2}\,
        &\lesssim\,
        \sum_{j\in\mathbb{N}}\,2^{(\nu-2\sigma)j}\,
        \int_{\mathbb{X}}\diff{x}\,\langle{x}\rangle^{-\nu}\,
        \Big|
        \int_{\mathbb{X}}\diff{y}\,
        \chi_{0}(\tfrac{|y^{-1}x|}{2^{j+2}})\,
        \chi_{\infty}(2|y^{-1}x|)\,
        k_{\sigma}(y^{-1}x)\,
        \widetilde{\chi}_{j}(y)\,f(y)
        \Big|^{2}\\[5pt]
        &=\,
        \sum_{j\in\mathbb{N}}\,2^{(\nu-2\sigma)j}\,
        \int_{\mathbb{X}}\diff{x}\,\langle{x}\rangle^{-\nu}\,
        \big|(\widetilde{\chi}_{j}\,f)\,*\,
        (\chi_{0}(\tfrac{|\cdot|}{2^{j+2}})\,
        \chi_{\infty}(2|\cdot|)\,k_{\sigma})
        \big|^{2}\\[5pt]
        &\lesssim\,
        \sum_{j\in\mathbb{N}}\,2^{(\nu-2\sigma)j}\,
        \|\chi_{0}(\tfrac{|\cdot|}{2^{j+2}})\,
        \chi_{\infty}(2|\cdot|)\,k_{\sigma}\|_{L^{2}(\mathbb{X})}^{2}\,
        \|\widetilde{\chi}_{j}f\|_{L^{2}(\mathbb{X})}^{2},
    \end{align*}
where
    \begin{align*}
        \|\chi_{0}(\tfrac{|\cdot|}{2^{j+2}})\,
        \chi_{\infty}(2|\cdot|)\,
        k_{\sigma}\|_{L^{2}(\mathbb{X})}^{2}\,
        &\asymp\,
        \int_{\mathbb{X}}\diff{x}\,
        \chi_{0}(\tfrac{|x|}{2^{j+2}})\,
        \chi_{\infty}(2|x|)\,
        |x|^{2\sigma-2\nu}\,\varphi_{0}^{2}(x)\\[5pt]
        &\lesssim\,
        \int_{1\le|x^{+}|\le2^{j}}\diff{x^{+}}\,
        |x^{+}|^{2\sigma-\nu-\ell}\\[5pt]
        &=\,\int_{1}^{2^{j}}\,\diff{r}\,r^{2\sigma-\nu-1}\,
        \lesssim\,2^{(2\sigma-\nu)j}.
    \end{align*}

We finally obtain
    \begin{align*}
         \|\mathcal{T}_{6}f(x)\|_{L^{2}(\mathbb{X})}^{2}\,
         \lesssim\,
         \sum_{j\in\mathbb{N}}\,
         \|\widetilde{\chi}_{j}f\|_{L^{2}(\mathbb{X})}^{2}\,
         \asymp\,\|f\|_{L^{2}(\mathbb{X})}^{2},
    \end{align*}
provided that $\gamma_{1},\gamma_{2}<\frac{\nu}{2}$ and $\gamma_{1}+\gamma_{2}\ge\sigma>0$. We omit the similar proof for the operator $\mathcal{T}_{7}$ and conclude that, for all $\gamma_{1},\gamma_{2}<\frac{\nu}{2}$ satisfying $\gamma_{1}+\gamma_{2}\ge\sigma>0$,
    \begin{align}\label{S3 estim T6T7}
        \|\mathcal{T}_{k}f\|_{L^{2}(\mathbb{X})}\,
        \lesssim\,\|f\|_{L^{2}(\mathbb{X})}
        \qquad\textnormal{with}\,\,\,k=6,\,7.
    \end{align}

\noindent\textbf{Conclusion:}
Since $\bigcup_{1\le{k}\le7}\supp\psi_{k}$ covers $\mathbb{X}$, we deduce, from \eqref{S3 estim T1}, \eqref{S3 estim T2T3}, \eqref{S3 estim T4T5}, and \eqref{S3 estim T6T7} that
    \begin{align*}
        \|\mathcal{T}f\|_{L^{2}(\mathbb{X})}\,
        \lesssim\,\sum_{1\le{k}\le7}
            \|\mathcal{T}_{k}f\|_{L^{2}(\mathbb{X})}\,
        \lesssim\,\|f\|_{L^{2}(\mathbb{X})}
    \end{align*}
    provided that $\sigma>0$, $\gamma_{1},\gamma_{2}<\min\lbrace{\frac{n}{2},\frac{\nu}{2}}\rbrace$, and $\sigma=\gamma_{1}+\gamma_{2}$.
\end{proof}

\begin{remark}
If we assume that $f$ is in addition bi-$K$-invariant in Theorem \ref{main thm SW} and Proposition \ref{S3 SW prop}, the last two cases in the above proof will be simplified according to the following trick:
    \begin{align}
        \int_{G}\diff{y}\,\varphi_{0}(y^{-1}x)\,f(y)
        =\,\int_{G}\diff{y}\,f(y)\,
            \int_{K}\diff{k}\,\varphi_{0}(y^{-1}kx)\,
        =\,\varphi_{0}(x)\,
            \int_{G}\diff{y}\,\varphi_{0}(y)\,f(y),
        \label{S32 SplitSph}
    \end{align}
since $\diff{k}$ is a normalized measure on $K$ and $\varphi_{0}$ is a spherical function. In fact, \eqref{S32 SplitSph} implies that
\begin{align*}
    \mathcal{T}_{j}(x)\,
    &=\,
    \int_{\mathbb{X}}\diff{x}\,
    \psi_{k}(x,y)\,|x|^{-\gamma_{1}}\,k_{\sigma}(y^{-1}x)\,
            |y|^{-\gamma_{2}}\,f(y)\\[5pt]
    &\asymp\,\varphi_{0}(x)\,
    \int_{\mathbb{X}}\diff{x}\,\psi_{k}(x,y)\,
        |x|^{-\gamma_{1}}\,|y^{-1}x|^{\sigma-\nu}\,
        |y|^{-\gamma_{2}}\,\varphi_{0}(y)\,f(y)
\end{align*}
for each $4\le{k}\le7$. Then we can conclude by using Lemma \ref{S3 SW lemma}.
\end{remark}

\begin{remark}\label{S3 rmk SWoptimal}
The regularity conditions $\gamma_{1},\gamma_{2}<\min\lbrace{\frac{n}{2},\frac{\nu}{2}}\rbrace$ in Theorem \ref{main thm SW} are necessary. In the Euclidean setting, the necessity of conditions occurring in the Stein-Weiss inequality is well explained in the recent note \cite{Ngo21}. On symmetric spaces, the kernel $k_{\sigma}$ behaves differently depending on whether it is close to or away from the origin. Hence the manifold dimension and the pseudo-dimension both play essential roles. To check this, it is sufficient to show that, for any  $\gamma_{1},\gamma_{2}\ge\min\lbrace{\frac{n}{2},\frac{\nu}{2}}\rbrace$, the double integral
    \begin{align}
        \int_{\mathbb{X}}\diff{x}\,\Big|
        \int_{\mathbb{X}}\diff{y}\,
            |x|^{-\gamma_{1}}\,k_{\sigma}(y^{-1}x)\,
            |y|^{-\gamma_{2}}\,f(y)
        \Big|^{2}
    \label{S32 doubleInt}
    \end{align}
is not finite for $x$ or $y$ located in some specific regions. Let us define a subset of $\mathfrak{a}^{+}$ consisting of vectors away from the walls:
    \begin{align*}
        \mathfrak{a}_{1}\,
        =\,\lbrace{
        H\in\mathfrak{a}^{+}\,|\,
        \langle{\alpha,H}\rangle\asymp|H|,\,\,\,\forall\,\alpha\in\Sigma^{+}
        }\rbrace.
    \end{align*}
Then, for any vector $H\in\mathfrak{a}_{1}$, we have 
\begin{align*}
    \delta(H)\,
    \gtrsim\,
    \begin{cases}
        |H|^{n-\ell}
        \qquad&\textnormal{if $|H|$ is bounded from above},\\[5pt]
        e^{\langle{2\rho,H}\rangle}
        \qquad&\textnormal{if $|H|$ is bounded from below}.
    \end{cases}
\end{align*}

\begin{itemize}[leftmargin=*]
    \item
    Let $f$ be a cut-off function such that $\supp{f}=\lbrace{y\in\mathbb{X}\,|\,\frac14\le|y|\le\frac12}\rbrace$. Notice that, for all $0<\sigma<n$, $|x|\le\frac12$, and $\frac14\le|y|\le\frac12$, we have
        \begin{align*}
            k_{\sigma}(y^{-1}x)\,
            \asymp\,
            |y^{-1}x|^{\sigma-n}\,
            \gtrsim\,
            |y|^{\sigma-n}
        \end{align*}
    according to \eqref{S3 Riesz ker estim}. Hence,
        \begin{align*}
            \eqref{S32 doubleInt}\,
            &\gtrsim\,
            \int_{\lbrace{x\,\in\,K(\exp{\mathfrak{a}_{1}})K\,|\,|x|\le\frac12}\rbrace}
            \diff{x}\,|x|^{-2\gamma_{1}}\,
            \Big|\int_{K(\exp{\mathfrak{a}_{1}})K\,\cap\,\supp{f}}\diff{y}\,
            |y|^{\sigma-n-\gamma_{2}}\Big|^{2}\\[5pt]
            &\gtrsim\,
            \int_{|x^{+}|\le\frac12}\diff{x^{+}}\,|x^{+}|^{-2\gamma_{1}+n-\ell}\,
            \Big|
            \int_{\frac14\le|y^{+}|\le\frac12}\diff{y^{+}}\,
            |y^{+}|^{\sigma-\gamma_{2}-\ell}
            \Big|^{2}\\[5pt]
            &=\,
            \Big\lbrace{
            \int_{0}^{\frac12}\diff{r}\,r^{-2\gamma_{1}+n-1}
            }\Big\rbrace\,
            \Big\lbrace{\underbrace{
            \int_{\frac14}^{\frac12}\diff{r}\,r^{\sigma-\gamma_{2}-1}}_{=\,\const}
            }\Big\rbrace,
        \end{align*}
    where the first integral is not finite for any $\gamma_{1}\ge\frac{n}{2}$. The necessity of $\gamma_{2}<\frac{n}{2}$ can be handled in the same way.
    
    \item
    Let $f$ be a bi-$K$-invariant cut-off function such that $\supp{f}=\lbrace{y\in\mathbb{X}\,|\,\frac{1}{2}\le|y|\le1}\rbrace$. For all $0<\sigma<\nu$, $|x|\ge2$, and $\frac{1}{2}\le|y|\le1$, we have $1\le|y^{-1}x|\le2|x|$ and then
        \begin{align*}
            k_{\sigma}(y^{-1}x)\,
            \asymp\,
            |y^{-1}x|^{\sigma-\nu}\,\varphi_{0}(y^{-1}x)\,
            \gtrsim\,
            |x|^{\sigma-\nu}\,\varphi_{0}(y^{-1}x)\,
        \end{align*}
    according to \eqref{S3 Riesz ker estim} again. Since $f$ is bi-$K$-invariant, we deduce from \eqref{S32 SplitSph} and \eqref{S2 phi0} that
        \begin{align*}
        \eqref{S32 doubleInt}\,
        &\gtrsim\,
            \int_{\lbrace{x\,\in\,K(\exp{\mathfrak{a}_{1}})K\,|\,|x|\ge2}\rbrace}
            \diff{x}\,|x|^{-2\gamma_{1}+2\sigma-2\nu}\,
            \Big|
            \int_{\frac12\le|y|\le1}\diff{y}\,
            |y|^{-\gamma_{2}}\,\varphi_{0}(y^{-1}x)\,
            \Big|^{2}\\[5pt]
        &=\,
            \int_{\lbrace{x\,\in\,K(\exp{\mathfrak{a}_{1}})K\,|\,|x|\ge2}\rbrace}
            \diff{x}\,|x|^{2\gamma_{2}-2\nu}\,\varphi_{0}^{2}(x)\,
            \underbrace{\Big|
            \int_{\frac12\le|y|\le1}\diff{y}\,
            |y|^{-\gamma_{2}}\,\varphi_{0}(y)
            \Big|^{2}}_{=\,\const}\\[5pt]
        &\gtrsim\,
            \int_{|x^{+}|\ge2}\diff{x^{+}}\,
            |x^{+}|^{2\gamma_{2}-\nu-\ell}\,
        \end{align*}
    which is not finite for any $\gamma_{2}\ge\frac{\nu}{2}$. Similarly, we can show that $\gamma_{1}<\frac{\nu}{2}$ is a necessary condition as well.
\end{itemize}
\end{remark}

\begin{remark}\label{S3 RmkSWCor}
In the proof of Proposition \ref{S3 SW prop}, the condition $\gamma_{1}+\gamma_{2}=\sigma$ is only required in the first two cases. In other two cases, it is sufficient to conclude with $\gamma_{1}+\gamma_{2}\ge\sigma$. If we consider the inequality with \textit{inhomogeneous} weights instead of the homogeneous ones, we can get rid of the analysis around the origin, then establish an inhomogeneous version of Stein-Weiss inequality with relaxed conditions $\gamma_{1},\gamma_{2}<\frac{\nu}{2}$ and $\gamma_{1}+\gamma_{2}\ge\sigma$, see the following Corollary. This type of inequality has been considered in \cite{Kai14} for $\gamma_{1}+\gamma_{2}>\sigma$. The endpoint improvement 
here allows one to deduce straightforwardly the Type (III) smoothing estimate in the full regularity range, see Theorem \ref{main thm SchSmoothing} and \cite[Remark 4.1]{Kai14}.
\end{remark}

\begin{corollary}\label{S3 SWCor}
Let $\sigma>0$, $\gamma_{1},\gamma_{2}<\frac{\nu}{2}$, and $\gamma_{1}+\gamma_{2}\ge\sigma$. Then we have
    \begin{align}\label{S3 inSW}
        \|\langle{\cdot}\rangle^{-\gamma_{1}}\,D_{x}^{-\sigma}\,
            f\|_{L^{2}(\mathbb{X})}\,
        \lesssim\,
        \|\langle{\cdot}\rangle^{\gamma_{2}}\,f\|_{L^{2}(\mathbb{X})}.
    \end{align}
\end{corollary}

\begin{proof}
Let us show that the operator $\widetilde{\mathcal{T}}:L^{2}(\mathbb{X})\rightarrow{L}^{2}(\mathbb{X})$ defined by
    \begin{align*}
    \widetilde{\mathcal{T}}f(x)\,=\,
        \int_{\mathbb{X}}\diff{y}\,
            \langle{x}\rangle^{-\gamma_{1}}\,k_{\sigma}(y^{-1}x)\,
            \langle{y}\rangle^{-\gamma_{2}}\,f(y)
    \end{align*}
is $L^{2}$-bounded. On the one hand, if $|x|$ and $|y|$ are both large, then $\langle{x}\rangle^{-\gamma_{1}}\langle{y}\rangle^{-\gamma_{2}}\asymp|x|^{-\gamma_{1}}|y|^{-\gamma_{2}}$ for any $\gamma_{1},\gamma_{2}\in\mathbb{R}$ and we go back to cases (i) and (iv) in the proof of Proposition \ref{S3 SW prop}. Notice that in Case (i), if one considers the inhomogeneous weight $\langle{x}\rangle$ instead of the homogeneous one, it is sufficient to take $j\in\mathbb{N}$ instead of $j\in\mathbb{Z}$ in the dyadic decomposition \eqref{S3 dyadic}, then conclude with relaxed condition $\gamma_{1}+\gamma_{2}\ge\sigma>0$, see \eqref{S3 estim T1} and \eqref{S3 estim T6T7}. On the other hand, if $|x|$ and $|y|$ are both small, then $|y^{-1}x|$ is also small and $\langle{x}\rangle^{-\gamma_{1}}\langle{y}\rangle^{-\gamma_{2}}$ is bounded for any $\gamma_{1},\gamma_{2}\in\mathbb{R}$. By using the the Kunze-Stein phenomenon \eqref{S3 KS1} and the kernel estimate \eqref{S3 Riesz ker estim}, we obtain
    \begin{align}\label{S3 Jap small1}
        \int_{\mathbb{X}}\diff{x}\,
        \Big|\int_{\mathbb{X}}\diff{y}\,
            \chi_{0}(|y^{-1}x|)\,k_{\sigma}(y^{-1}x)\,f(y)\Big|^{2}\,
        &=\,\|(\chi_{0}k_{\sigma})*f\|_{L^{2}(\mathbb{X})}^{2}\notag\\[5pt]
        &\lesssim\,\|\chi_{0}k_{\sigma}\varphi_{0}\|_{L^{1}(\mathbb{X})}^{2}\, 
            \|f\|_{L^{2}(\mathbb{X})}^{2}
    \end{align}
where
    \begin{align}\label{S3 Jap small2}
        \|\chi_{0}k_{\sigma}\varphi_{0}\|_{L^{1}(\mathbb{X})}\, 
        \lesssim\,
            \int_{|x^{+}|\le1}\diff{x}\,
            |x^{+}|^{\sigma-n}\,|x^{+}|^{n-\ell}\,
        <\,+\infty
    \end{align}
for any $\sigma>0$. 

In the remaining cases where $|x|$ and $|y|$ are not both small or both large, we go back to Case (iii) in the proof of Proposition \ref{S3 SW prop}. Notice that, in contrast to $|x|^{-\gamma_{1}}$, the inhomogeneous weight $\langle{x}\rangle^{-\gamma_{1}}$ has no contribution when $|x|$ is small. Hence the condition $\gamma_{1}<\frac{n}{2}$ is not required. Similarly, we can remove the condition $\gamma_{2}<\frac{n}{2}$ as well when $|y|$ is small. 

We deduce that $\widetilde{\mathcal{T}}$ is $L^{2}$-bounded, provided that $\gamma_{1},\gamma_{2}<\frac{\nu}{2}$ and $\gamma_{1}+\gamma_{2}\ge\sigma>0$.
\end{proof}

%%%%%%%%%%%%%%%%%%%%%%%%%%%%%%%%%%%%%%%%%%%%%%%%%%%%%%%%%%%%%%%%%%%%%%%%%%%%%%%%%%
\subsection{Resolvent estimate and smoothing property}
We prove in this part the resolvent estimate stated in Theorem \ref{main thm resolv}. Combining it with the Stein-Weiss inequality Theorem \ref{main thm SW}, we deduce Theorem \ref{main thm smoothing}. As we mentioned in the introduction, the standard scaling argument carried out in \cite{KaYa89,Sug03} fails in the present setting. We prove Theorem \ref{main thm resolv} along the lines in \cite{Chi08,Kai14} with a more careful analysis around the origin, since we are considering the homogeneous weights. We will combine the improved $L^2$-continuity of $\mathcal{JR}$-transform (Lemma \ref{S2 JR continuity}) with two estimates borrowing from the Euclidean Fourier analysis. Recall that $\mathfrak{a}$ is an $\ell$-dimensional flat submanifold of $\mathbb{X}$. Let $g$ be a reasonable function on $\mathfrak{a}$. Then we have the following.
\begin{itemize}
    \item (Besov embedding)
        If $\ell=1$ and $\frac12<\theta<\frac32$, we have
        \begin{align}\label{S3 Morrey}
        |g(\lambda_{1})-g(\lambda_{2})|\,
        \lesssim\,|\lambda_{1}-\lambda_{2}|^{\theta-\frac12}\,
            \||\cdot|^{\theta}\,
                \mathcal{F}_{\mathfrak{a}}g\|_{L^{2}(\mathfrak{a})}
            \qquad\forall\,\lambda_{1},\lambda_{2}\in\mathfrak{a}.
        \end{align}
        When $\theta=1$, the estimate \eqref{S3 Morrey} is the classical Morrey inequality. For other $\frac12<\theta<\frac32$, it is sufficient to notice that the Hölder-Zygmund space with index $0<\theta-\frac12<1$ is the standard Besov space $B_{\infty\infty}^{\theta-1/2}(\mathfrak{a})$, which is embedded into $B_{22}^{\theta}(\mathfrak{a})$. See, for instance, \cite[Sect. 2.2.2]{Saw18}.
        
    \vspace{5pt}\item (Fourier restriction theorem)
        If $\ell\ge2$ and $\frac12<\theta<\frac{\ell}{2}$, we have
        \begin{align}\label{S3 Restriction}
        \int_{|\lambda|=\,r}\diff{\sigma}_{\lambda}\,
        |(\mathcal{F}_{\mathfrak{a}}g)(\lambda)|^{2}\,
        \,\lesssim\,
        r^{2\theta-1}\,\int_{\mathfrak{a}}\diff{\lambda}\,
        |\lambda|^{2\theta}\,|g(\lambda)|^{2}
        \qquad\forall\,r>0,
    \end{align}
see for instance \cite[Theorem 5.6]{BlSa92}. Here $\diff{\sigma}_{\lambda}$ denotes the usual surface measure.
\end{itemize}

\begin{proof}[Proof of Theorem \ref{main thm resolv}]
According to the Plancherel formula and the transform \eqref{S2 transformXtoR}, we write
    \begin{align*}
        (D^{2\alpha}(D^{2}-\zeta)^{-1}\,f,\,f)_{L^{2}(\mathbb{X})}\,
        =\,|W|^{-1}\,\int_{B}\diff{b}\,\int_{\mathfrak{a}}\diff{\lambda}\,
            |\mathbf{c}(\lambda)|^{-2}\,|\mathcal{F}f(\lambda,b)|^{2}\,
            \frac{|\lambda|^{2\alpha}}{|\lambda|^{2}-\zeta},
    \end{align*}
where
    \begin{align*}
        |\mathbf{c}(\lambda)|^{-2}\,|\mathcal{F}f(\lambda,b)|^{2}\,
        =\,\big||\mathbf{c}(\lambda)|^{-1}\,
                \mathcal{F}_{\mathfrak{a}}
                [\mathcal{R}f(\cdot,b)](\lambda)\big|^{2}\,
        =\,\big|\mathcal{F}_{\mathfrak{a}}
                [\mathcal{JR}f(\cdot,b)](\lambda)\big|^{2},
    \end{align*}
for all $\lambda\in\mathfrak{a}$ and $b\in{B}$. 

When $\mathbb{X}$ is of rank $\ell=1$, by applying \eqref{S3 Morrey} with $g(\lambda)=\mathbf{c}(\lambda)^{-1}\mathcal{F}f(\lambda,B)$, $\lambda_{2}=0$ and $\theta=1-\alpha$, we have
    \begin{align*}
        |\mathbf{c}(\lambda)|^{-2}\,|\mathcal{F}f(\lambda,B)|^{2}\,
        \lesssim\,
        |\lambda|^{1-2\alpha}\,
            \||\cdot|^{1-\alpha}\,
                \mathcal{JR}f(\cdot,B)\|_{L^{2}(\mathfrak{a})}^{2},
    \end{align*}
provided that $-\frac12<\alpha<\frac12$. When $\mathbb{X}$ is of rank $\ell\ge2$, we use \eqref{S3 Restriction} with $\theta=1-\alpha$, and obtain
    \begin{align*}
    (D^{2\alpha}(D^{2}-\zeta)^{-1}\,f,\,f)_{L^{2}(\mathbb{X})}\,
    &=\,|W|^{-1}\,\int_{B}\diff{b}\,\int_{0}^{+\infty}\diff{r}\,
            \frac{r^{2\alpha}}{r^{2}-\zeta}\,
            \int_{|\lambda|=\,r}\diff{\sigma}_{\lambda}\,
            \big|\mathcal{F}_{\mathfrak{a}}
                [\mathcal{JR}f(\cdot,b)](\lambda)\big|^{2}\\[5pt]
    &\lesssim\,|W|^{-1}\,\int_{B}\diff{b}\,\int_{0}^{+\infty}\diff{r}\,
            \frac{r}{r^{2}-\zeta}\,
            \int_{\mathfrak{a}}\diff{\lambda}\,
            |\lambda|^{2-2\alpha}\,
            \big|\mathcal{JR}f(\lambda,b)\big|^{2}
    \end{align*}
provided that $1-\frac{\ell}{2}<\alpha<\frac{1}{2}$. Hence, for all $\ell\ge1$, we have
    \begin{align*}
        |\im(D^{2\alpha}(D^{2}-\zeta)^{-1}\,f,\,f)_{L^{2}(\mathbb{X})}|\,
        &\lesssim\,
            \|\mathcal{JR}f\|_{L^{2}(\mathfrak{a}\times{B},
            \,|\lambda|^{2-2\alpha}\diff{\lambda}\diff{b})}^{2}
            \Big|\im\int_{0}^{+\infty}\diff{r}\,
                \frac{r}{r^{2}-\zeta}\Big|\\[5pt]
        &\,\lesssim\,\|f\|_{L^{2}(\mathbb{X},|x|^{2-2\alpha}\diff{x})},
    \end{align*}
where 
    \begin{align*}
        \Big|\im\int_{0}^{+\infty}\diff{r}\,\frac{r}{r^{2}-\zeta}\Big|\,
        \le\,\frac{1}{2}\,\int_{0}^{+\infty}\diff{s}\,
            \frac{|\im\zeta|}{(s-\re\zeta)^{2}+(\im\zeta)^{2}}\,
        \le\,\frac{\pi}{2}\,
        \qquad\forall\,\zeta\in\mathbb{C}\smallsetminus\mathbb{R}.
    \end{align*}
We conclude that, for suitable $\alpha$, 
    \begin{align*}
        |\im(D^{2\alpha}(D^{2}-\zeta)^{-1}\,f,\,f)_{L^{2}(\mathbb{X})}|\,
        \lesssim\,\|f\|_{L^{2}(\mathbb{X},|x|^{2-2\alpha}\diff{x})},
    \end{align*}
according to Lemma \ref{S2 JR continuity}.
\end{proof}

\begin{remark}
We give here only an estimate for the imaginary part. According to Remark \ref{S1 Kato}, this is enough to deduce the smoothness for the corresponding operator, namely, the smoothing property for the free Schrödinger equation. The real part estimate can be handled along the lines in \cite{Chi08,Kai14}, using the Stein-Weiss inequality \eqref{main thm SW ineq} instead of theirs. This allows one to prove the so-called super-smoothness for the corresponding operator, see \cite{Kat66,KaYa89}.
\end{remark}

Combining the resolvent estimate \eqref{main thm resolv} with the Stein-Weiss inequality \eqref{main thm SW ineq}, we deduce our Theorem \ref{main thm smoothing}.

\begin{proof}[Proof of Theorem \ref{main thm smoothing}]
According to Remark \ref{S1 Kato} and Theorem \ref{main thm resolv}, we know that the smoothing property \eqref{main thm smoothing schrodinger} holds for all $-\frac12<\alpha<\frac12$ in rank one ($\ell=1$) and $1-\frac{\ell}{2}<\alpha<\frac12$ in higher ranks ($\ell\ge2$).
It remains for us to show that the smoothing property \eqref{main thm smoothing schrodinger} remains valid for all $1-\min\lbrace{\frac{n}{2},\frac{\nu}{2}}\rbrace<\alpha\le1-\frac{\ell}{2}$ when $\ell\ge2$. Let $0<\varepsilon<\frac12$ be a small constant. We write
    \begin{align*}
        \||x|^{\alpha-1}\,D_{x}^{\alpha}\,u
        \|_{L^{2}(\mathbb{R}_{t}\times\mathbb{X})}\,
        =\,
        \||x|^{-(1-\alpha)}\,D_{x}^{-(\frac12-\varepsilon-\alpha)}\,
        D_{x}^{\frac12-\varepsilon}\,u
        \|_{L^{2}(\mathbb{R}_{t}\times\mathbb{X})}
    \end{align*}
where $\frac12-\varepsilon-\alpha\ge\frac{\ell}{2}-\frac{1}{2}-\varepsilon>\frac{\ell}{2}-1\ge0$ and $1-\alpha<\min\lbrace{\frac{n}{2},\frac{\nu}{2}}\rbrace$ fulfil conditions of the Stein-Weiss inequality \eqref{main thm SW ineq}. Hence, we obtain
    \begin{align*}
        \||x|^{\alpha-1}\,D_{x}^{\alpha}\,u
        \|_{L^{2}(\mathbb{R}_{t}\times\mathbb{X})}\,
        \lesssim\,
        \||x|^{-\frac12-\varepsilon}\,D_{x}^{\frac12-\varepsilon}\,
        u\|_{L^{2}(\mathbb{R}_{t}\times\mathbb{X})}.
    \end{align*}
Notice that
\begin{align*}
    \||x|^{-\frac12-\varepsilon}\,D_{x}^{\frac12-\varepsilon}\,
        u\|_{L^{2}(\mathbb{R}_{t}\times\mathbb{X})}\,
    =\,
    \||x|^{(\frac12-\varepsilon)-1}\,D_{x}^{\frac12-\varepsilon}\,
        u\|_{L^{2}(\mathbb{R}_{t}\times\mathbb{X})},
\end{align*}
where $1-\frac{\ell}{2}\le0<\frac{1}{2}-\varepsilon<\frac12$ for all $\ell\ge2$. We deduce that the smoothing property
    \begin{align*}
        \||x|^{\alpha-1}\,D_{x}^{\alpha}\,u
        \|_{L^{2}(\mathbb{R}_{t}\times\mathbb{X})}\,
        \lesssim\,
        \|u_{0}\|_{L^{2}(\mathbb{X})}
    \end{align*}
remains valid for all $1-\min\lbrace{\frac{n}{2},\frac{\nu}{2}}\rbrace<\alpha\le1-\frac{\ell}{2}$ when $\ell\ge2$. We conclude that the smoothing property \eqref{main thm smoothing schrodinger} holds for all $1-\min\lbrace{\frac{n}{2},\frac{\nu}{2}}\rbrace<\alpha<\frac12$ in higher ranks.
\end{proof}

\begin{remark}\label{S3 remark optimal}
On symmetric spaces $G/K$ where $G$ is complex, the regularity condition in Theorem \ref{main thm smoothing} is optimal, i.e., the smoothing property 
    \begin{align*}
        \||x|^{\alpha-1}\,D_{x}^{\alpha}\,u
        \|_{L^{2}(\mathbb{R}_{t}\times\mathbb{X})}\,
        \lesssim\,
        \|u_{0}\|_{L^{2}(\mathbb{X})}
    \end{align*}
cannot hold for any $\alpha\le1-\frac{\nu}{2}$ or $\alpha\ge\frac12$. When $G$ is complex, we can write the spherical Fourier transform \eqref{S22 HarishChandra} as
    \begin{align*}
        f(x)\,=\,\const\,\varphi_{0}(x)\,
        \int_{\mathfrak{p}}\diff{\lambda}\,
        \mathcal{H}f(\lambda)\,e^{-i\langle{\lambda,x}\rangle}
        \qquad\forall\,f\in\mathcal{S}(K\backslash{G/K}),
    \end{align*}
where $\mathfrak{p}$ is an $n$-dimensional flat space, see \cite[Theorem 4.7 and Theorem 9.1]{Hel00}. Let $u_{0}$ be a bi-$K$-invariant function such that its spherical Fourier transform is radial in $\mathfrak{p}$. Then 
    \begin{align*}
        D_{x}^{\alpha}\,e^{itD_{x}^{2}}u_{0}(x)\,
        &=\,\const\,\varphi_{0}(x)\,
        \int_{\mathfrak{p}}\diff{\lambda}\,
        |\lambda|^{\alpha}\,e^{-it|\lambda|^{2}}\,
        (\mathcal{H}u_{0})(|\lambda|)\,
        e^{-i\langle{\lambda,x}\rangle}\\[5pt]
        &=\,\const\,\varphi_{0}(x)\,
        \int_{0}^{\infty}\diff{r}\,
        r^{\alpha}\,e^{-itr^{2}}\,
        (\mathcal{H}u_{0})(r)\,
        \int_{|\lambda|=r}\diff{\sigma_{\lambda}}\,
        e^{-i\langle{\lambda,x}\rangle},
    \end{align*}
where the inner integral is the modified Bessel function:
\begin{align*}
    \int_{|\lambda|=r}\diff{\sigma_{\lambda}}\,
        e^{-i\langle{\lambda,x}\rangle}\,
    =\,r^{\frac{n}{2}}\,|x|^{\frac{2-n}{2}}\,J_{\frac{n-2}{2}}(r|x|).
\end{align*}
By making the change of variable $r=\sqrt{s}$, we obtain
    \begin{align*}
        D_{x}^{\alpha}\,e^{itD_{x}^{2}}u_{0}(x)\,
        =\,\const\,|x|^{\frac{2-n}{2}}\,\varphi_{0}(x)\,
        \int_{0}^{\infty}\diff{s}\,e^{-its}
        s^{\frac{\alpha}{2}+\frac{n}{4}-\frac{1}{2}}\,
        (\mathcal{H}u_{0})(\sqrt{s})\,J_{\frac{n-2}{2}}(\sqrt{s}|x|).
    \end{align*}
Together with the Plancherel formula (in variable $t$), we deduce that
    \begin{align}
        &\||x|^{\alpha-1}\,D_{x}^{\alpha}\,e^{itD^{2}}u_{0}
            \|_{L^2(\mathbb{R}_{t}\times\mathbb{X})}^{2}
            \notag\\[5pt]
        &=\,\const\,
        \int_{\mathbb{X}}\diff{x}\,|x|^{2\alpha-n}\,\varphi_{0}^{2}(x)\,
        \int_{0}^{\infty}\diff{s}\,s^{\alpha+\frac{n}{2}-1}\,
        |(\mathcal{H}u_{0})(\sqrt{s})|^{2}\,
        \Big|J_{\frac{n-2}{2}}(\sqrt{s}|x|)\Big|^{2}
        \notag\\[5pt]
        &=\,\const\,
        \int_{0}^{\infty}\diff{s}\,s^{\alpha+\frac{n}{2}-1}\,
        |(\mathcal{H}u_{0})(\sqrt{s})|^{2}\,
        \int_{\mathbb{X}}\diff{x}\,|x|^{2\alpha-n}\,\varphi_{0}^{2}(x)\,
        \Big|J_{\frac{n-2}{2}}(\sqrt{s}|x|)\Big|^{2}.
        \label{S31 innerintegral}
    \end{align}
For any fixed $s>0$, the Bessel function $|J_{\frac{n-2}{2}}(\sqrt{s}|x|)|$ behaves asymptotically as $|x|^{\frac{n-2}{2}}$ when $|x|\rightarrow{0}$ and $|x|^{-{1}/{2}}$ when $|x|\rightarrow\infty$. Moreover, notice on the one hand that
    \begin{align*}
        \int_{\lbrace{x\,\in\,{K(\exp\mathfrak{a}_{1})K}\,|\,|x|\le1}\rbrace}\diff{x}\,
        |x|^{2\alpha-n}\,\varphi_{0}^{2}(x)\,|x|^{n-2}\,
        &\gtrsim\,
        \int_{|x^{+}|\le1}\diff{x^{+}}\,|x^{+}|^{2\alpha-2}\,|x|^{n-\ell}\\[5pt]
        &=\,
        \int_{0}^{1}\diff{r}\,r^{2\alpha-2+n-1}
    \end{align*}
which is finite if and only if $\alpha>1-\frac{n}{2}$. On the other hand,
    \begin{align*}
        \int_{\lbrace{x\,\in\,{K(\exp\mathfrak{a}_{1})K}\,|\,|x|\ge1}\rbrace}
        \diff{x}\,|x|^{2\alpha-n}\,\varphi_{0}^{2}(x)\,|x|^{-1}\,
        &\gtrsim\,
        \int_{|x^{+}|\ge1}\diff{x^{+}}\,|x^{+}|^{2\alpha-n+\nu-\ell-1}\\[5pt]
        &=\,
        \int_{1}^{+\infty}\diff{r}\,r^{2\alpha-n+\nu-2}
    \end{align*}
which is finite provided that $\alpha<\frac{\nu-n}{2}+\frac{1}{2}$. Here $\mathfrak{a}_{1}\subset\mathfrak{a}^{+}$ is the subset consisting of vectors away from the walls, see Remark \ref{S3 rmk SWoptimal}. Since $n=\nu$ in the case where $G$ is complex, then the inner integral on the right hand side of \eqref{S31 innerintegral} is finite if and only if $1-\frac{\nu}{2}<\alpha<\frac{1}{2}$. Hence, when $G$ is complex, the smoothing property \eqref{main thm smoothing schrodinger} cannot hold for any $\alpha\le1-\frac{\nu}{2}$ or $\alpha\ge\frac12$.
\end{remark}

%%%%%%%%%%%%%%%%%%%%%%%%%%%%%%%%%%%%%%%%%%%%%%%%%%%%%%%%%%%%%%%%%%%%%%%%%%%%%%%%%%
%%%%%%%%%%                            SECTION IV                      %%%%%%%%%%%%
%%%%%%%%%%%%%%%%%%%%%%%%%%%%%%%%%%%%%%%%%%%%%%%%%%%%%%%%%%%%%%%%%%%%%%%%%%%%%%%%%%
\section{Comparison principle on symmetric spaces}\label{Section.4 Comparison}
Consider two evolution equations corresponding to operators $a_{1}(D_{x})$ and $a_{2}(D_{x})$:
    \begin{align*}
        \begin{cases}
            (i\partial_{t}+a_{1}(D_{x}))\,u(t,x)=\,0,\\[5pt]
            u(0,x)\,=\,u_{0}(x),
        \end{cases}
        \qquad\textnormal{and}\qquad
        \begin{cases}
            (i\partial_{t}+a_{2}(D_{x}))\,v(t,x)=\,0,\\[5pt]
            v(0,x)\,=\,v_{0}(x),
        \end{cases}
    \end{align*}
whose solutions are given by $u(t,x)=e^{ita_{1}(D_{x})}u_{0}(x)$ and $v(t,x)=e^{ita_{2}(D_{x})}v_{0}(x)$. The comparison principle allows one to compare smoothing properties between these two different equations when the symbols of $a_{1}(D_{x})$ and $a_{2}(D_{x})$ satisfy certain relations. We extend this tool to symmetric spaces along the lines in \cite[Theorem 2.5]{RuSu12}, since most of our arguments are made in the Cartan subspace $\mathfrak{a}$, which is an $\ell$-dimensional flat submanifold of $\mathbb{X}$.

\begin{theorem}[First comparison principle]
\label{S4 Comparison principle}
Let $\tau_{1},\tau_{2}$ be two continuous functions on $\mathbb{R}_{+}$. Let $a_{1},a_{2}\in\mathcal{C}^{1}(\mathbb{R}_{+})$ be real-valued and strictly monotone functions on the support of a measurable function $\chi$ on $\mathbb{R}_{+}$. If there exist some constants $C>0$ such that two pairs of functions $\lbrace{\tau_{1},a_{1}}\rbrace$ and $\lbrace{\tau_{2},a_{2}}\rbrace$ fulfil the comparison condition
    \begin{align}\label{S4 CC}
        \frac{|\tau_{1}(r)|}{|a_{1}'(r)|^{1/2}}\,
        \le\,C\,\frac{|\tau_{2}(r)|}{|a_{2}'(r)|^{1/2}}
        \tag{CC}
    \end{align}
for all $r\in\supp\chi$ satisfying $a_{1}'(r)\neq0$ and $a_{2}'(r)\neq0$, then for any measurable function $\omega$ on $\mathbb{X}$, we have
    \begin{align}\label{S4 CC weight}
        \|\omega(x)\,\chi(D_{x})\,\tau_{1}(D_{x})\,e^{ita_{1}(D_{x})}\,
        &u_{0}(x)\|_{L^{2}(\mathbb{R}_{t}\times\mathbb{X})}\notag\\[5pt]
        &\le\,C\,
        \|\omega(x)\,\chi(D_{x})\,\tau_{2}(D_{x})\,e^{ita_{2}(D_{x})}\,
        u_{0}(x)\|_{L^{2}(\mathbb{R}_{t}\times\mathbb{X})},
    \end{align}
where the equality holds if \eqref{S4 CC} holds with equality.
\end{theorem}

\begin{remark}
When functions $a_1$ and $a_2$ satisfy \eqref{S4 CC} for all $r\in\mathbb{R}_{+}$, Theorem \ref{S4 Comparison principle} and Corollary \ref{S4 Secondary comparison principle} hold globally without function $\chi$. The reason to introduce such a function into the estimates is that the comparison relation between symbols may vary in different frequencies. This is the case for wave-type equations, see Corollary \ref{corollary freq}.
\end{remark}

\begin{proof}
As usual, we assume that all the integrals below make sense, we perform calculation on the set $a_{1}'(r)\neq0$, where the inverse of $a_{1}$ is differentiable. By using the inverse formula \eqref{S2 Inverse Helgason Fourier} of the Helgason-Fourier transform and polar coordinates, we write
    \begin{align*}
        &\chi(D_{x})\,\tau_{1}(D_{x})\,
        e^{ita_{1}(D_{x})}\,u_{0}(x)\\[5pt]
        &=\,
        |W|^{-1}\,\int_{\mathfrak{a}}\diff{\lambda}\,|\mathbf{c}(\lambda)|^{-2}\,
            (\chi\tau_{1})(|\lambda|)\,e^{ita_{1}(|\lambda|)}\,
            \int_{B}\diff{b}\,e^{\langle{i\lambda+\rho,\,A(k^{-1}x)}\rangle}\,
            \mathcal{F}u_{0}(\lambda,b)\\[5pt]
        &=\,
        |W|^{-1}\,\int_{0}^{+\infty}\diff{r}\,
            r^{\ell-1}\,(\chi\tau_{1})(r)\,e^{ita_{1}(r)}\,
            \underbrace{
            \int_{\mathbb{S}^{\ell-1}}\,\diff{\sigma_{\eta}}\,|\mathbf{c}(r\eta)|^{-2}\,
            \int_{B}\,\diff{b}\,
            e^{\langle{ir\eta+\rho,\,A(k^{-1}x)}\rangle}\,\mathcal{F}u_{0}(r\eta,b)
            }_{=\,U(r,x)}.
    \end{align*}
By substituting $r=a_{1}^{-1}(s)$ on the support of $\chi$, we have
    \begin{align*}
        &\chi(D_{x})\,\tau_{1}(D_{x})\,
        e^{ita_{1}(D_{x})}\,u_{0}(x)\\[5pt]
        &=\,
        |W|^{-1}\,\int_{a_{1}(\mathbb{R}_{+})}\diff{s}\,
            |(a_{1}^{-1})'(s)|\,(a_{1}^{-1}(s))^{\ell-1}\,
            (\chi\tau_{1})(a_{1}^{-1}(s))\,
            U(a_{1}^{-1}(s),x)\,e^{its}.
    \end{align*}
Applying the Plancherel formula (in variable $t$), we obtain
    \begin{align*}
        &\|\chi(D_{x})\,\tau_{1}(D_{x})\,e^{ita_{1}(D_{x})}\,
            u_{0}(x)\|_{L^{2}(\mathbb{R})}^{2}\\[5pt]
        &=\,(2\pi)^{-1}|W|^{-2}\,
            \int_{a_{1}(\mathbb{R}_{+})}\,\diff{s}\,
            |(a_{1}^{-1})'(s)|^{2}\,|a_{1}^{-1}(s)|^{2\ell-2}\,
            |(\chi\tau_{1})(a_{1}^{-1}(s))|^{2}\,
            |U(a_{1}^{-1}(s),x)|^{2}\\[5pt]
        &=\,(2\pi)^{-1}|W|^{-2}\,
            \int_{0}^{+\infty}\,\diff{r}\,
                r^{2\ell-2}\,|\chi(r)|^{2}\,
                \frac{|\tau_{1}(r)|^{2}}{|a_{1}'(r)|}\,
                |U(r,x)|^{2}.
    \end{align*}
Here, we have used the substitution $s=a_{1}(r)$ and the identity $(a_{1}^{-1})'(a_{1}(r))=a_{1}'(r)^{-1}$. We deduce from the comparison condition \eqref{S4 CC} that
    \begin{align*}
        &\|\chi(D_{x})\,\tau_{1}(D_{x})\,
        e^{ita_{1}(D_{x})}\,u_{0}(x)\|_{L^{2}(\mathbb{R}_{t})}\\[5pt]
        &\le\,
        C\,(2\pi)^{-1}|W|^{-2}\,
            \int_{0}^{+\infty}\,\diff{r}\,
                r^{2\ell-2}\,|\chi(r)|^{2}\,
                \frac{|\tau_{2}(r)|^{2}}{|a_{2}'(r)|}\,
                |U(r,x)|^{2}\\[5pt]
        &=\,C\,
        \|\chi(D_{x})\,\tau_{2}(D_{x})\,
        e^{ita_{2}(D_{x})}\,u_{0}(x)\|_{L^{2}(\mathbb{R}_{t})}
    \end{align*}
for all $x\in\mathbb{X}$. Then $\eqref{S4 CC weight}$ follows.
\end{proof}

%%%%%%%%%%%%%%%%%%%%%%%%%%%%%%%%%%%%%%%%%%%%%%%%%%%%%%%%%%%%%%%%%%%%%%%%%%%%%%
\subsection{General order Schrödinger-type equations}
The above comparison principle allows us to deduce some new smoothing estimates from the model case. Consider the Schrödinger-type equation of order $m>0$: 
    \begin{align}\label{S4 SchroOrder}
        (i\partial_{t}+D_{x}^{m})\,u(t,x)=\,0,
        \qquad\,u(0,x)\,=\,u_{0}(x),
    \end{align}    
whose solution is given by $u(t,x)=e^{itD_{x}^{m}}u_{0}(x)$. We will show that the solution to the Cauchy problem \eqref{S4 SchroOrder} satisfies the smoothing property:
    \begin{align}\label{S4 SmoothSchA}
        \|A(x,D_{x})\,u\|_{L^{2}(\mathbb{R}_{t}\times\mathbb{X})}\,
        \lesssim\,
        \|u_{0}\|_{L^{2}(\mathbb{X})},
    \end{align}
namely, Theorem \ref{main thm SchSmoothing}. Recall that $A(x,D)$ is defined as one of the following:

\begin{table}[ht]
\setlength{\tabcolsep}{5pt}
\renewcommand{\arraystretch}{2}
\begin{tabular}{|c|c|c|c|}
\hline
\cellcolor{gray!25} Type
& \cellcolor{gray!25} $A(x,D)$ 
& \cellcolor{gray!25} $\ell=1$ 
& \cellcolor{gray!25} $\ell\ge2$ \\
\hline
\textnormal{(I)}
& $|x|^{\alpha-\frac{m}{2}}D^{\alpha}$
& $\frac{m-3}{2}<\alpha<\frac{m-1}{2}$
& $\frac{m-\min\lbrace{n,\nu}\rbrace}{2}<\alpha<\frac{m-1}{2}$\\ 
\hline
\textnormal{(II)}
&$\langle{x}\rangle^{-s}D^{\frac{m-1}{2}}$
& \multicolumn{2}{c|}{$s>\frac{1}{2}$ \textnormal{and} $m>0$} \\  
\hline
\textnormal{(III)}
&$\langle{x}\rangle^{-s}\langle{D}\rangle^{\frac{m-1}{2}}$
& \multicolumn{2}{c|}{{$s\ge\frac{m}{2}$ \textnormal{and} $1<m<\nu$}} \\
\hline
\end{tabular}
\vspace{10pt}
\caption{Regularity conditions on $\mathbb{X}$ for the Schrödinger-type equations with order $m>0$.}
\label{S4 TableSch}
\end{table}

Let us clarify what is already known about Theorem \ref{main thm SchSmoothing} and what remains for us to prove. The Type (I) estimate is proved for $m=2$ in Theorem \ref{main thm smoothing}. In \cite{Kai14}, the author showed
\begin{itemize}[itemsep=5pt]
    \item the Type (II) estimate for $m>1$;
    \item the Type (III) estimate for $\ell\ge2$ and $1<m<\ell$;
    \item the Type (III) estimate for $1<m\le\nu$, but with restricted condition $s>\frac{m}{2}$.
\end{itemize}
Notice that if the Type (II) estimate holds for $m=1$, then the Type (III) estimate holds for $s>\frac{m}{2}$ with $m=1$. As we mentioned in Remark \ref{S3 RmkSWCor}, the Type (III) estimate in the critical (optimal) case $s=\frac{m}{2}$ is a consequence of the improved Stein-Weiss inequality \eqref{S3 inSW}, see also \cite[Remark 4.1]{Kai14}. 

It remains for us to prove the Type (I) estimate for all $m>0$ and complete the Type (II) estimate when $0<m\le1$. They are the consequences of the comparison principle.

\begin{proof}[Proof of Theorem \ref{main thm SchSmoothing}]
Let $-\frac{1}{2}<\alpha'<\frac12$ when $\ell=1$ and $1-{\min\lbrace{\frac{n}{2},\frac{\nu}{2}}\rbrace}<\alpha'<\frac{1}{2}$ when $\ell\ge2$. Notice that, for any $m>0$ and $r>0$, the two pairs of functions
\begin{align*}
    \lbrace{\tau_{1}(r)=r^{\frac{m}{2}+\alpha'-1},a_{1}(r)=r^{m}}\rbrace
    \qquad\textnormal{and}\qquad
    \lbrace{\tau_{2}(r)=r^{\alpha'},a_{2}(r)=r^{2}}\rbrace
\end{align*}
satisfy the comparison condition \eqref{S4 CC} with $C=\sqrt{\frac{2}{m}}$. Hence
    \begin{align*}
        \||x|^{\alpha'-1}\,D_{x}^{\frac{m}{2}+\alpha'-1}\,
            &e^{-itD_{x}^{m}}\,u_{0}(x)\|_{
        L^{2}(\mathbb{R}_{t}\times\mathbb{X})}\\[5pt]
        &=\,\sqrt{\frac{2}{m}}\,
        \||x|^{\alpha'-1}\,D_{x}^{\alpha'}\,
            e^{-itD_{x}^{2}}\,u_{0}(x)\|_{
        L^{2}(\mathbb{R}_{t}\times\mathbb{X})}\,
        \lesssim\,
        \|u_{0}\|_{L^{2}(\mathbb{X})}
    \end{align*}
according to Theorem \ref{S4 Comparison principle} and Theorem \ref{main thm smoothing}. We deduce the Type (I) smoothing estimate by taking $\alpha=\frac{m}{2}+\alpha'-1$.

Now, suppose that $0<m\le1<m'$. For any $r>0$, the two pairs of functions  $\lbrace\tau_{1}(r)=r^{{(m-1)}/{2}},$ $a_{1}(r)=r^{m}\rbrace$ and $\lbrace{\tau_{2}(r)=r^{{(m'-1)}/{2}},a_{2}(r)=r^{m'}}\rbrace$ satisfy \eqref{S4 CC} with $C=\sqrt{\frac{m}{m'}}$. Then, we obtain
    \begin{align*}
        \|\langle{x}\rangle^{-s}\,|D_{x}|^{\frac{m-1}{2}}\,
            &e^{-itD_{x}^{m}}\,u_{0}(x)\|_{
        L^{2}(\mathbb{R}_{t}\times\mathbb{X})}\\[5pt]
        &=\,\sqrt{\frac{m}{m'}}\,
        \|\langle{x}\rangle^{-s}\,|D_{x}|^{\frac{m'-1}{2}}\,
            e^{-itD_{x}^{m'}}\,u_{0}(x)\|_{
        L^{2}(\mathbb{R}_{t}\times\mathbb{X})}\,
        \lesssim\,
        \|u_{0}\|_{L^{2}(\mathbb{X})}
    \end{align*}
for all $s>\frac12$. Hence, the Type (II) smoothing estimate holds for all $m>0$.
\end{proof}

\subsection{Wave and Klein-Gordon equations}
Another important example is the Cauchy problem
    \begin{align}\label{S4 KGequ}
        \begin{cases}
            (\partial_{t}^{2}+D_{x}^{2}+\zeta)\,u(t,x)=\,0,\\[5pt]
            u(0,x)\,=\,u_{0}(x),\,\partial_{t}|_{t=0}\,u(t,x)\,=\,u_{1 }(x),
        \end{cases}
    \end{align}
which is the wave equation when $\zeta=0$ and the Klein-Gordon equation when $\zeta>0$. To establish the smoothing properties for \eqref{S4 KGequ}, we introduce the following secondary comparison principle. 

\begin{corollary}[Secondary comparison principle]
\label{S4 Secondary comparison principle}
Suppose that $s>1/2$ and $\alpha$ satisfies
    \begin{align}\label{S4 RankCdt}
        \begin{cases}
        -\frac12<\alpha<\frac12
        &\qquad\textnormal{if $\ell=1$},\\[5pt]
        1-\frac{\min\lbrace{n,\nu}\rbrace}{2}<\alpha<\frac12
        &\qquad\textnormal{if $\ell\ge2$}. 
        \end{cases}
    \end{align}
Let $a\in\mathcal{C}^{1}(\mathbb{R}_{+})$ be a real-valued and strictly monotone function on the support of a measurable function $\chi$ on $\mathbb{R}_{+}$. Let $\tau\in\mathcal{C}^{0}(\mathbb{R}_{+})$ be such that, for some $C>0$, we have
    \begin{align}\label{S4 SCC}
        |\tau(r)|\,\le\,C\,\sqrt{|a'(r)|}
        \tag{SCC}
    \end{align}
for all $r\in\supp\chi$. Then, the solution to the Cauchy problem
\begin{align}\label{S4 SecEqu}
    (i\partial_{t}+a(D_{x}))u(t,x)\,=\,0,
    \qquad\,u(t,x)\,=\,u_{0}(x),
\end{align}
satisfies
    \begin{align}
        \|\langle{x}\rangle^{-s}\,\chi(D_{x})
            \tau(D_{x})\,e^{ita(D_{x})}\,
        u_{0}(x)\|_{L^{2}(\mathbb{R}_{t}\times\mathbb{X})}\,
        &\lesssim\,\|u_{0}\|_{L^{2}(\mathbb{X})},
        \label{S4 SCC1}\\[5pt]
        \||x|^{\alpha-1}\,\chi(D_{x})\,D_{x}^{\alpha-\frac{1}{2}}\,
        \tau(D_{x})\,e^{ita(D_{x})}\,
        u_{0}(x)\|_{L^{2}(\mathbb{R}_{t}\times\mathbb{X})}\,
        &\lesssim\,\|u_{0}\|_{L^{2}(\mathbb{X})}.
        \label{S4 SCC2}
    \end{align} 
\end{corollary}

\begin{proof}
This corollary is a straightforward consequence of the first comparison principle and the smoothing property of the Schrödinger equation. Notice that if $\tau(r)$ and $a(r)$ satisfy the second comparison condition \eqref{S4 SCC} with $C=\frac{1}{\sqrt{2}}$, then $\lbrace{\tau_{1}(r)=r^{\alpha-1/2}\tau(r),a_{1}(r)=a(r)}\rbrace$ and $\lbrace{\tau_{2}(r)=r^{\alpha},a_{2}(r)=r^{2}}\rbrace$ fulfill \eqref{S4 CC} with the same constant. Hence, for all $\alpha$ satisfies \eqref{S4 RankCdt}, we have
    \begin{align*}
        \||x|^{\alpha-1}\,\chi(D_{x})\,D_{x}^{\alpha-\frac{1}{2}}\,
            &\tau(D_{x})\,
            e^{ita(D_{x})}\,u_{0}(x)\|_{L^{2}(\mathbb{R}_{t}\times\mathbb{X})}\\[5pt]
        &\lesssim\,
        \||x|^{\alpha-1}\,\chi(D_{x})\,D_{x}^{\alpha}\,
            e^{itD_{x}^{2}}\,u_{0}(x)\|_{L^{2}(\mathbb{R}_{t}\times\mathbb{X})}\,
        \lesssim\,\|u_{0}\|_{L^{2}(\mathbb{X})}
    \end{align*}
according to Theorem \ref{S4 Comparison principle} and Theorem \ref{main thm smoothing}. 

Similarly, since $\lbrace{\tau(r),a(r)}\rbrace$ and $\lbrace{\tau_{3}(r)=r^{{(m-1)}/{2}},a_{3}(r)=r^{m}}\rbrace$ satisfy \eqref{S4 CC} for any $m>0$, we deduce from Theorem \ref{S4 Comparison principle} and Theorem \ref{main thm smoothing} again that 
    \begin{align*}
        \|\langle{x}\rangle^{-s}\,\chi(D_{x})
            &\tau(D_{x})\,e^{ita(D_{x})}\,
        u_{0}(x)\|_{L^{2}(\mathbb{R}_{t}\times\mathbb{X})}\\[5pt]
        &\lesssim\,
        \|\langle{x}\rangle^{-s}\,\chi(D_{x})
            D_{x}^{\frac{m-1}{2}}\,e^{itD_{x}^{m}}\,
        f(x)\|_{L^{2}(\mathbb{R}_{t}\times\mathbb{X})}\,
        \lesssim\,\|u_{0}\|_{L^{2}(\mathbb{X})}
    \end{align*}
for any $s>\frac12$.
\end{proof}

By applying the above corollary with $a(r)=\sqrt{r^{2}+\mu}$, where $\mu\ge0$, we obtain the following smoothing estimates in different frequencies. Notice that 
    \begin{align*}
        \sqrt{|a'(r)|}\,\ge\,
        \begin{cases}
            1
            &\quad\,\textnormal{if}\,\,\,\mu=0\,\,\,
            \textnormal{or}\,\,\,r>1,\\[5pt]
            \sqrt{r}
            &\quad\,\textnormal{if}\,\,\,0<r\le1.
        \end{cases}
    \end{align*}

\begin{corollary}\label{corollary freq}
Suppose that $s>1/2$ and $\alpha$ satisfies \eqref{S4 RankCdt}. Let $\chi$ be a smooth cut-off function on $\mathbb{R}_{+}$ such that $\chi=1$ around the origin and denote by $U_{l}=\chi(D)u_{0}$ and $U_{h}=(1-\chi(D))u_{0}$ the initial data corresponding to the low and high frequencies. Then, for all $\mu\ge0$,
    \begin{align}
        \||x|^{\alpha-1}D^{\alpha}\,
        e^{\pm{it}\sqrt{D_{x}^{2}+\mu}}\,U_{l}(x)\|_{
        L^{2}(\mathbb{R}_{t}\times\mathbb{X})}\,
        &\lesssim\,
        \|U_{l}\|_{L^{2}(\mathbb{X})},
        \label{S4 low1}\\[5pt]
        \||x|^{\alpha-1}D^{\alpha-\frac{1}{2}}\,
        e^{\pm{it}\sqrt{D_{x}^{2}+\mu}}\,U_{h}(x)\|_{
        L^{2}(\mathbb{R}_{t}\times\mathbb{X})}\,
        &\lesssim\,
        \|U_{h}\|_{L^{2}(\mathbb{X})},
        \label{S4 high1}\\[5pt]
        \|\langle{x}\rangle^{-s}\,D_{x}^{\frac{1}{2}}\,
        e^{\pm{it}\sqrt{D_{x}^{2}+\mu}}\,U_{l}(x)\|_{
        L^{2}(\mathbb{R}_{t}\times\mathbb{X})}\,
        &\lesssim\,
        \|U_{l}\|_{L^{2}(\mathbb{X})},
        \label{S4 low2}\\[5pt]
        \|\langle{x}\rangle^{-s}\,
        e^{\pm{it}\sqrt{D_{x}^{2}+\mu}}\,U_{h}(x)\|_{
        L^{2}(\mathbb{R}_{t}\times\mathbb{X})}\,
        &\lesssim\,
        \|U_{h}\|_{L^{2}(\mathbb{X})}.
        \label{S4 high2}
    \end{align}
Moreover, in the limiting case where $\mu=0$, we have better estimates in the low-frequency part:
    \begin{align}
        \||x|^{\alpha-1}D_{x}^{\alpha-\frac{1}{2}}\,
        e^{\pm{it}D_{x}}\,U_{l}(x)\|_{
        L^{2}(\mathbb{R}_{t}\times\mathbb{X})}\,
        &\lesssim\,
        \|U_{l}\|_{L^{2}(\mathbb{X})},
        \label{S4 low3}\\[5pt]
        \|\langle{x}\rangle^{-s}\,
        e^{\pm{it}D_{x}}\,U_{l}(x)\|_{
        L^{2}(\mathbb{R}_{t}\times\mathbb{X})}\,
        &\lesssim\,
        \|U_{l}\|_{L^{2}(\mathbb{X})}.
        \label{S4 low4}
    \end{align}
\end{corollary}

The usual way to relate smoothing estimates of wave and Schrödinger equations relies on changing variables in the corresponding restriction theorems. The previous corollary allows us to relate them simply according to the comparison principle. By combining estimates \eqref{S4 high1} and \eqref{S4 low3}, as while as \eqref{S4 high2} and \eqref{S4 low4}, we deduce the following smoothing estimates for the wave equation.

\begin{theorem}\label{S4 wave thm}
Consider the Cauchy problem \eqref{S4 KGequ} with $\zeta=0$, namely, the wave equation. We have the smoothing properties
    \begin{align}
        \|\langle{x}\rangle^{-s}\,u\|_{
        L^{2}(\mathbb{R}_{t}\times\mathbb{X})}\,
        &\lesssim\,
        \|u_{0}\|_{L^{2}(\mathbb{X})}\,+\,\|D_{x}^{-1}\,u_{1}\|_{L^{2}(\mathbb{X})},
        \label{S4 wave1}\\[5pt]
        \||x|^{\beta-\frac{1}{2}}\,D_{x}^{\beta}\,u\|_{
        L^{2}(\mathbb{R}_{t}\times\mathbb{X})}\,
        &\lesssim\,
        \|u_{0}\|_{L^{2}(\mathbb{X})}\,+\,\|D_{x}^{-1}\,u_{1}\|_{L^{2}(\mathbb{X})},
        \label{S4 wave2}
    \end{align}
for any $s>\frac12$ and $\beta$ satisfies
    \begin{align}
        \begin{cases}
            -1<\beta<0
            &\qquad\textnormal{if}\,\,\,\ell=1, \\[5pt]
            \frac{1-\min\lbrace{n,\nu}\rbrace}{2}<\beta<0
            &\qquad\textnormal{if}\,\,\,\ell\ge2.
        \end{cases}
    \label{S4 RankCdtbeta}
    \end{align}
\end{theorem}

Similar estimates such as \eqref{S4 wave1} and \eqref{S4 wave2} are well-known in $\mathbb{R}^{N}$ for $N\ge2$ and $\frac{1-N}{2}<\beta<0$, see \cite{Ben94,RuSu12}. In particular, since the regularity range of the Kato-type smoothing property of the Schrödinger equation is wider on $\mathbb{H}^2$ (see Remark \ref{S1 Kato}), analogous phenomenon appears in studying the wave equation: estimate \eqref{S4 wave2} holds for all $-1<\beta<0$ on $\mathbb{H}^2$, while similar estimate holds on $\mathbb{R}^2$ if and only if $-\frac{1}{2}<\beta<0$.

By combining \eqref{S4 high2} and \eqref{S4 low1} (with $\alpha=0$), we can deduce the following smoothing property for the Klein-Gordon equation. 

\begin{theorem}\label{S4 KG thm}
Consider the Cauchy problem \eqref{S4 KGequ} with $\zeta>0$, namely, the Klein-Gordon equation. We have the smoothing property
    \begin{align}\label{S4 smoothingKG}
        \|\langle{x}\rangle^{-1}\,u\|_{
        L^{2}(\mathbb{R}_{t}\times\mathbb{X})}\,
        &\lesssim\,
        \|u_{0}\|_{L^{2}(\mathbb{X})}\,+\,\|D_{x}^{-1}\,u_{1}\|_{L^{2}(\mathbb{X})}.
    \end{align}
\end{theorem}

Similar estimate as \eqref{S4 smoothingKG} has been established in $\mathbb{R}^{N}$ for $N\ge3$, see \cite{Ben94,RuSu12}. Moreover, if one considers the weight $\langle{x}\rangle^{-s}$ with $s>1$ instead of $\langle{x}\rangle^{-1}$, similar smoothing property remains valid on $\mathbb{R}^{2}$. In our setting, the estimate \eqref{S4 smoothingKG} always holds even in the $2$-dimensional case. The reason is that $|x|^{-1}$ is $D^{2}$-smooth on $\mathbb{H}^{2}$, which is not the case on $\mathbb{R}^{2}$.

%%%%%%%%%%%%%%%%%%%%%%%%%%%%%%%%%%%%%%%%%%%%%%%%%%%%%%%%%%%%%%%%%%%%%%%%%%%%%%%%
\subsection{Other examples}\label{subsection other examples}
Many other equations in the Euclidean setting admit the smoothing properties, but are less considered on more general manifolds because of the lack of physical backgrounds. From the point of view of mathematical analysis, the above arguments also allow us to deduce their smoothing properties on symmetric spaces easily. The following are two examples.

\begin{corollary}[Smoothing estimate of the relativistic Schrödinger equation]\label{S4 relaSch}
Consider the Cauchy problem
\begin{align*}
    \begin{cases}
        (i\partial_{t}-\sqrt{1+D_{x}^{2}})\,u(t,x)=\,0,\\[5pt]
        u(0,x)\,=\,u_{0}(x).
    \end{cases}
\end{align*}
Then, we have the smoothing property 
    \begin{align}
        \|\langle{x}\rangle^{-1}\,u\|_{
        L^{2}(\mathbb{R}_{t}\times\mathbb{X})}\,
        \lesssim\,
        \|u_{0}\|_{L^{2}(\mathbb{X})}.
        \label{S4 rela}
    \end{align}
\end{corollary}

Analogous estimate holds in $\mathbb{R}^{N}$ for all $N\ge3$, and the order of its weight $\langle{x}\rangle^{-1}$ is sharp, see \cite{BeNe97,Wal02}. Notice that we do not need the limiting absorption principle used in \cite{BeNe97}. The smoothing property \eqref{S4 rela} is a straightforward consequence of estimates \eqref{S4 high2} and \eqref{S4 low1}. Therefore, we have similar phenomena to the Klein-Gordon equation on $\mathbb{H}^{2}$.

Beyond the Schrödinger-type equation with \textit{constant coefficients}, we can also deduce similar smoothing properties for some \textit{time-variable coefficients} equations on symmetric spaces. Consider the Cauchy problem
    \begin{align}\label{S4 SchroOrderTime}
        (i\partial_{t}+\theta'(t)D_{x}^{m})\,u(t,x)=\,0,
        \qquad\,u(0,x)\,=\,u_{0}(x),
    \end{align}
where $\theta$ is a suitable function on $\mathbb{R}$. In \cite{FeRu20}, the authors established the comparison principle and some smoothing properties for  \eqref{S4 SchroOrderTime} in $\mathbb{R}^{N}$ when $\theta$ satisfies 
\begin{align}\label{S4 timecoefhyp}
    \begin{cases}
        \theta\in\mathcal{C}^{1}(\mathbb{R}),\\
        \theta(0)=0,\\
        \textnormal{$\theta$ is strictly monotone or $\theta'$ vanishes at finitely many points}.
    \end{cases}
\end{align}
See also \cite{KPRV05,MMT12,CiRe14,FeSt21} and the references therein for equations with other variable coefficients. Notice that the Cauchy problem \eqref{S4 SchroOrderTime} with $\theta$ satisfying \eqref{S4 timecoefhyp} covers the Schrödinger equation \eqref{S4 SchroOrder} if one sets $\theta(t)=t$. The following analogous property is a straightforward consequence of \cite[Lemma 2.1]{FeRu20} and Theorem \ref{main thm SchSmoothing}.

\begin{corollary}\label{S4 timeSch}
Suppose that $\theta$ meets the condition \eqref{S4 timecoefhyp} and $A(x,D)$ is described as in Table \ref{S1 TableSch}. Then, we have
    \begin{align}
        \||\theta'(t)|^{\frac{1}{2}}\,A(x,D_{x})\,e^{i\theta(t)D_{x}^{m}}\,u_{0}\|_{
        L^{2}(\mathbb{R}_{t}\times\mathbb{X})}\,
        \lesssim\,
        \|u_{0}\|_{L^{2}(\mathbb{X})}.
    \end{align}
   
\end{corollary}

%%%%%%%%%%%%%%%%%%%%%%%%%%%%%%%%%%%%%%%%%%%%%%%%%%%%%%%%%%%%%%%%%%%%%%%%%%%%%%%%%%
%%%%%%%%%%                             APPENDIX                       %%%%%%%%%%%%
%%%%%%%%%%%%%%%%%%%%%%%%%%%%%%%%%%%%%%%%%%%%%%%%%%%%%%%%%%%%%%%%%%%%%%%%%%%%%%%%%%
\appendix
\section{Proofs of two lemmas}
We prove Lemma \ref{S2 JR continuity} and Lemma \ref{S3 SW lemma} in this appendix.

\begin{lemma}
For any $\sigma\ge0$, we have
    \begin{align}
        \|\mathcal{JR}f\|_{L^{2}(\mathfrak{a}\times{B},\,
            |H|^{2\sigma}\diff{H}\diff{b})}\,
        \lesssim\,
        \|f\|_{L^{2}(\mathbb{X},\,|x|^{2\sigma}\diff{x})}.
    \end{align}
\end{lemma}
%%%%%
\begin{proof}
Since we consider here the homogeneous weights, we modify slightly the cut-off functions that we used in the proof of Lemma \ref{S3 KS lemma}. Let $\chi_{0}\in\mathcal{C}_{c}^{\infty}(\mathbb{R}_{+})$ be a cut-off function such that $\chi_{0}=1$ on $[0,\frac34]$ and $\chi_{0}=0$ on $[1,+\infty]$. For every $j\in\mathbb{Z}$, we define
    \begin{align*}
        \phi_{j}(r)\,
        =\,\chi_{0}(2^{-j}r)\,-\,\chi_{0}(2^{-j+1}r)
        \qquad\forall\,r\ge0,
    \end{align*}
which is compactly supported in $[3\times2^{j-3},2^{j}]$ and satisfies $\phi_{j}=1$ on $[2^{j-1},3\times2^{j-2}]$. Moreover, we have
    \begin{align}\label{S2 JR Dyadic}
        \sum_{j\in\mathbb{Z}}\,\phi_{j}(r)\,=\,1
        \qquad\forall\,r\ge0.
    \end{align}
Notice that
    \begin{align*}
        \|f\|_{L^{2}(\mathbb{X},\,|x|^{2\sigma}\diff{x})}^{2}\,
        \asymp\,
            \sum_{j\in\mathbb{Z}}\,2^{2\sigma{j}}\,
            \|\phi_{j}(|\cdot|)\,f\|_{L^{2}(\mathbb{X})}^{2}
    \end{align*}
and
    \begin{align*}    
        \|\mathcal{JR}f\|_{L^{2}(\mathfrak{a}\times{B},\,
            |H|^{2\sigma}\diff{H}\diff{b})}^{2}\,
        \asymp\,
            \sum_{j\in\mathbb{Z}}\,2^{2\sigma{j}}\,
            \|\phi_{j}(|\cdot|)\,
            \mathcal{JR}f\|_{L^{2}(\mathfrak{a}\times{B})}^{2}.
    \end{align*}

Recall the factorization $|\mathbf{c}(\lambda)|^{-1}=\mathbf{b}(\lambda)^{-1}\bm\pi(i\lambda)$ where $\bm{\pi}(i\lambda)=\prod_{\alpha\in\Sigma_{r}^{+}}\langle{\alpha,\lambda}\rangle$ and $\mathbf{b}(\lambda)$ satisfies
    \begin{align}\label{S2 JR b0}
        |p(\tfrac{\partial}{\partial\lambda})\mathbf{b}(-\lambda)^{\pm1}|\,
        \lesssim\,
            \underbrace{\prod\nolimits_{\alpha\in\Sigma_{r}^{+}}\,
            (1+|\langle{\alpha,\lambda}\rangle|)^{
            \mp\frac{m_{\alpha}+m_{2\alpha}}{2}\pm1}
            }_{=\,\mathbf{b}_{0}(\lambda)^{\pm1}}
        \qquad\forall\lambda\in\mathfrak{a}+i\overline{\mathfrak{a}^{+}},
    \end{align}
for any differential polynomial $p(\tfrac{\partial}{\partial\lambda})$, see for instance \cite[pp.1041-1042]{AnJi99}. Using the partition of unity \eqref{S2 JR Dyadic}, we write, for all $H\in\mathfrak{a}$ and $b\in{B}$,
    \begin{align*}
        \mathcal{JR}f(H,b)\,
        &=\,\sum_{h\in\mathbb{Z}}\,
            |\mathbf{c}(D_{H})|^{-1}\,\mathcal{R}(\phi_{h}f)(H,b)\\[5pt]
        &=\,\sum_{h\in\mathbb{Z}}\,
            \mathbf{b}(D_{H})^{-1}\,\widetilde{\chi}_{h}(H)\,
            \pi(iD_{H})\mathcal{R}(\phi_{h}f)(H,b)\\[5pt]
        &=\,\sum_{h\in\mathbb{Z}}\,
            \mathbf{b}(D_{H})^{-1}\,\widetilde{\chi}_{h}(H)\,\mathbf{b}(D_{H})\,
            \mathcal{JR}(\phi_{h}f)(H,b)
    \end{align*}
where $\widetilde{\chi}_{h}$ is a cut-off function on $\mathfrak{a}$ such that $\widetilde{\chi}_{h}=1$ in the support of $\bm\pi(iD_{H})\mathcal{R}(\phi_{h}f)\subset\supp\phi_{h}$. Then, for all $H\in\supp\widetilde{\chi}_{h}$, we have $|H|\le|K(\exp{H})N|\lesssim2^{h}$. We deduce, for every $j\in\mathbb{Z}$, that
    \begin{align}\label{S2 JR estimate}
        &2^{\sigma{j}}\,\|\phi_{j}(|\cdot|)\,
            \mathcal{JR}f\|_{L^{2}(\mathfrak{a}\times{B})}\notag\\[5pt]
        &\lesssim\,
            2^{\sigma{j}}\,\sum_{h\in\mathbb{Z}}\,\|\phi_{j}(|\cdot|)\,
            \mathbf{b}^{-1}(D_{H})\,\widetilde{\chi}_{h}\,\mathbf{b}(D_{H})\,
             \mathcal{JR}(\phi_{h}f)\|_{L^{2}(\mathfrak{a}\times{B})}\notag\\[5pt]
        &\le\,
            \sum_{h\in\mathbb{Z}}\,
            \underbrace{
            2^{\sigma(j-h)}\,
            \|\phi_{j}(|\cdot|)\,
            \mathbf{b}^{-1}(D_{H})\,\widetilde{\chi}_{h}\,\mathbf{b}(D_{H})
            \|_{L^{2}(\mathfrak{a})\rightarrow{L}^{2}(\mathfrak{a})}^{2}
            }_{=\,a_{jh}}
            \big(2^{\sigma{h}}\|\phi_{h}(|\cdot|)\,f\|_{L^{2}(\mathbb{X})}\big), 
    \end{align}
since $\mathcal{JR}:L^{2}(\mathbb{X})\longrightarrow{L^{2}(\mathfrak{a}\times{B},|W|^{-1}\diff{H}\diff{b})}$ is an isometry.
%%%%%
We claim that there exist two constants $C_{1},C_{2}>0$, independent of $j,h\in\mathbb{Z}$, such that
    \begin{align}\label{S2 JR uniform ajk}
        \sum_{j\in\mathbb{Z}}\,a_{jh}\,\le\,C_{1}
        \qquad\textnormal{and}\qquad
        \sum_{h\in\mathbb{Z}}\,a_{jh}\,\le\,C_{2}.
    \end{align}

We obtain from \eqref{S2 JR estimate} and \eqref{S2 JR uniform ajk} that
    \begin{align*}
        \|\mathcal{JR}f\|_{L^{2}(\mathfrak{a}\times{B},\,
            |H|^{2\sigma}\diff{H}\diff{b})}^{2}\,
        &\asymp\,\sum_{j\in\mathbb{Z}}\,2^{2\sigma{j}}\,
            \|\phi_{j}(|\cdot|)\,
            \mathcal{JR}f\|_{L^{2}(\mathfrak{a}\times{B})}^{2}\\[5pt]
        &\lesssim\,
            \sum_{j\in\mathbb{Z}}\,
            \Big\lbrace{
            \sum_{h\in\mathbb{Z}}\,a_{jh}\,
            \big(2^{\sigma{h}}\|\phi_{h}(|\cdot|)\,f\|_{L^{2}(\mathbb{X})}\big)
            }\Big\rbrace^{2}\\[5pt]
        &\le\,
            \sum_{j\in\mathbb{Z}}\,
            \Big\lbrace{\sum_{h\in\mathbb{Z}}\,a_{jh}}\Big\rbrace
            \Big\lbrace{
            \sum_{h\in\mathbb{Z}}\,a_{jh}\,
            \big(2^{2\sigma{h}}\|\phi_{h}(|\cdot|)\,f\|_{L^{2}(\mathbb{X})}^{2}\big)
            }\Big\rbrace\\[5pt]
        &\le\,
            C_{1}C_{2}\,\sum_{h\in\mathbb{Z}}\,
                2^{2\sigma{h}}\|\phi_{h}(|\cdot|)f
                \|_{L^{2}(\mathbb{X})}^{2}\\[5pt]
        &\lesssim\,
        \|f\|_{L^{2}(\mathbb{X},\,|x|^{2\sigma}\diff{x})}^{2}.
    \end{align*}

It remains us to check \eqref{S2 JR uniform ajk}. We recall the Weyl-H\"ormander pseudo-differential calculus used in \cite{Kai14}. Let us denote by $S(m,g)$ the Weyl-H\"ormander symbol class for a slowly varying metric $g$ and a $g$-continuous positive function $m$. A differential operator is uniformly bounded from $L^2$ to $L^{2}$ if its symbol belongs to $S(m,g)$ with $m$ bounded, see \cite[Ch.XVIII]{Hor94}. According to \eqref{S2 JR b0}, we know that $\mathbf{b}(\lambda)^{\pm1}\in{S}(\mathbf{b}_{0}(\lambda)^{\pm1},\langle{H}\rangle^{-2}|\!\diff{H}|^{2}+|\!\diff{\lambda}|^{2})$. Moreover, for every $j,h\in\mathbb{Z}$, the differential symbols $\phi_{j}(|H|)$, $\widetilde{\chi}_{h}(H)$, $2^{(\sigma+1)j}\phi_{j}(|H|)|H|^{-(\sigma+1)}$, and $2^{-(\sigma+1)h}\widetilde{\chi}_{h}(H)|H|^{\sigma+1}$ are all in the class $S(1,\langle{H}\rangle^{-2}|\!\diff{H}|^{2}+|\!\diff{\lambda}|^{2})$. Hence, the families of pseudo-differential operators 
    \begin{align*}
        \lbrace{
            \phi_{j}(|\cdot|)\,\mathbf{b}^{-1}(D_{H})\,
            \widetilde{\chi}_{h}\,\mathbf{b}(D_{H})
        }\rbrace_{j,h\in\mathbb{Z}}\\[5pt]
        \lbrace{
            2^{(\sigma+1)(j-h)}
            \phi_{j}(|\cdot|)\,\mathbf{b}^{-1}(D_{H})\,
            \widetilde{\chi}_{h}\,\mathbf{b}(D_{H})
        }\rbrace_{j,h\in\mathbb{Z}}
    \end{align*}
are uniformly bounded from $L^{2}(\mathfrak{a})$ to $L^{2}(\mathfrak{a})$. From where we deduce that
    \begin{align*}
        \sum_{j\le{h}}\,a_{jh}\,
        =\,\sum_{j\le{h}}\,
            2^{\sigma(j-h)}
            \underbrace{
            \|\phi_{j}(|\cdot|)\,
            \mathbf{b}^{-1}(D_{H})\,\widetilde{\chi}_{h}\,\mathbf{b}(D_{H})
            \|_{L^{2}(\mathfrak{a})\rightarrow{L}^{2}(\mathfrak{a})}^{2}
            }_{\le\,C}\,
        \le\,\frac{2^{\sigma}C}{2^{\sigma}-1}
    \end{align*}
and 
    \begin{align*}
        \sum_{j\ge{h}}\,a_{jh}\,
        &=\,\sum_{j\ge{h}}\,2^{h-j}\,
            \overbrace{
            2^{(\sigma+1)(j-h)}
            \|\phi_{j}(|\cdot|)\,
            \mathbf{b}^{-1}(D_{H})\,\widetilde{\chi}_{h}\,\mathbf{b}(D_{H})
            \|_{L^{2}(\mathfrak{a})\rightarrow{L}^{2}(\mathfrak{a})}^{2}
            }^{\le\,C}\\[5pt]
        &\le\,C\sum_{m\in\mathbb{N}}2^{-m}\,\le\,2C,
    \end{align*}
which prove \eqref{S2 JR uniform ajk}. Then the proof is complete.
\end{proof}

\begin{lemma}
Let $\mathcal{K}:\mathbb{R}_{+}\times\mathbb{R}_{+}\rightarrow\mathbb{R}_{+}$ be a homogeneous function of degree $-\ell$ such that
    \begin{align*}
        \int_{0}^{+\infty}\diff{s}\,s^{\tfrac{\ell}{2}-1}\,\mathcal{K}(1,s)\,<\,\infty.
    \end{align*}
Let $\kappa_{1}$ and $\kappa_{2}$ be two bi-$K$-invariant functions on $G$ and let $C_{1},C_{2}$ be two constants such that 
    \begin{align*}
        |\kappa_{j}(\exp{x}^{+})|\,\delta^{1/2}(x^{+})\,
        \le\,C_{j}
        \qquad\forall\,x\in{G},\,\,
        \forall\,j=1,\,2.
    \end{align*}
Then the operator $S:L^{2}(G)\rightarrow{L^{2}}(G)$ defined by
    \begin{align*}
        Sf(x)\,
        =\,\kappa_{1}(x)\,
            \int_{G}\diff{y}\,
            \mathcal{K}(|x|,|y|)\,
            \kappa_{2}(y)\,f(y)
    \end{align*}
is bounded.

\end{lemma}

\begin{proof}
We prove this lemma along the lines in \cite[Lemma 2.1]{StWe58}. According to the Cartan decomposition, there exists a constant $C>0$ such that
    \begin{align*}
        Sf(x)\,
        =\,C\,\kappa_{1}(\exp{x^{+}})\,
            \int_{K}\diff{k_{1}}\,\int_{K}\diff{k_{2}}\,
            \underbrace{
            \int_{\mathfrak{a}^{+}}\diff{y^{+}}\,\delta(y^{+})\,
            \mathcal{K}(|x^{+}|,|y^{+}|)\,
            \kappa_{2}(\exp{y^{+}})\,f(k_{1}(\exp{y^{+}})k_{2})
            }_{=\,I_{k_{1},k_{2}}(|x^{+}|)}.
    \end{align*}
By using the polar coordinates $x^{+}=R\xi$ and $y^{+}=r\eta$ with $R,r>0$ and $\xi,\eta\in\mathbb{S}^{\ell-1}$, then the substitution $r=sR$ and the homogeneity of $\mathcal{K}$, we write, for every $k_{1},k_{2}\in{K}$,
    \begin{align}\label{main thm SW jensen}
       I_{k_{1},k_{2}}(R)\,
        &=\,\int_{\mathbb{S}^{\ell-1}}\diff{\sigma_{\eta}}\,
            \int_{0}^{\infty}\diff{r}\,r^{\ell-1}\,
            \delta(r\eta)\,\mathcal{K}(R,r)\,
            \kappa_{2}(\exp(r\eta))\,f(k_{1}(\exp{(r\eta))k_{2}})\notag\\[5pt]
        &=\,\int_{\mathbb{S}^{\ell-1}}\diff{\sigma_{\eta}}
            \underbrace{
            \int_{0}^{\infty}\diff{s}\,s^{\ell-1}\,
                \delta(sR\eta)\,\mathcal{K}(1,s)\,
                \kappa_{2}(\exp(sR\eta))\,f(k_{1}(\exp{(sR\eta))k_{2}})
            }_{=\,\widetilde{I}_{\eta,k_{1},k_{2}}(R)}.
    \end{align}

By duality, there exists a function $h$ on $\mathbb{R}_{+}$ such that $\int_{0}^{\infty}\diff{R}R^{\ell-1}|h(R)|^{2}=1$ and
    \begin{align*}
        \Big\lbrace
        \int_{0}^{\infty}\diff{R}\,R^{\ell-1}\,|\widetilde{I}_{\eta,k_{1},k_{2}}(R)|^{2}
        \Big\rbrace^{1/2}
        =\,
        \int_{0}^{\infty}\diff{R}\,R^{\ell-1}\,h(R)\,\widetilde{I}_{\eta,k_{1},k_{2}}(R).
    \end{align*}
We obtain, by using successively the Fubini theorem, the Cauchy-Schwarz inequality and the substitution $r=sR$, that
    \begin{align*}
        &\Big\lbrace
            \int_{0}^{\infty}\diff{R}\,R^{\ell-1}\,
            |\widetilde{I}_{\eta,k_{1},k_{2}}(R)|^{2}
            \Big\rbrace^{1/2}
        \notag\\[5pt]
        &\le\,
            \int_{0}^{\infty}\diff{s}\,s^{\ell-1}\,\mathcal{K}(1,s)\,
            \int_{0}^{\infty}\diff{R}\,R^{\ell-1}\,|h(R)\,
            \delta(sR\eta)\,\kappa_{2}(\exp(sR\eta))\,f(k_{1}(\exp{(sR\eta))k_{2}})|
            \notag\\[5pt]
        &\le\,
            \Big\lbrace{
            \int_{0}^{\infty}\diff{r}\,r^{\ell-1}\,
            \delta(r\eta)\,|f(k_{1}(\exp{(r\eta))k_{2}})|^{2}}
            \Big\rbrace^{1/2}
            \notag\\[5pt]
            &\hspace{50pt}
                \underbrace{
                \int_{0}^{\infty}\diff{s}\,s^{\ell-1}\,
                s^{-\tfrac{\ell}{2}}\mathcal{K}(1,s)\,
                \Big\lbrace{
                \int_{0}^{\infty}\diff{R}\,R^{\ell-1}\,
                \delta(sR\eta)\,|\kappa_{2}(\exp(sR\eta))|^{2}\,|h(R)|^{2}}\Big\rbrace^{1/2}
                }_{=\,\bigO(1)}     
    \end{align*}
since $\delta(sR\eta)\,|\kappa_{2}(\exp(sR\eta))|^{2}$ is uniformly bounded by assumption. By applying the Jensen inequality to \eqref{main thm SW jensen}, we deduce from the above estimate that
    \begin{align}
     \int_{0}^{\infty}\diff{R}\,R^{\ell-1}\,|I_{k_{1},k_{2}}(R)|^{2}\,
    &\le\,|\omega_{\ell-1}|\,
        \int_{0}^{\infty}\diff{R}\,R^{\ell-1}\,
            \int_{\mathbb{S}^{\ell-1}}\,\diff{\sigma_{\eta}}\,
            |\widetilde{I}_{\eta,k_{1},k_{2}}(R)|^{2}\notag\\[5pt]
    &\lesssim\,
            \int_{\mathbb{S}^{\ell-1}}\,\diff{\sigma_{\eta}}\,
            \int_{0}^{\infty}\diff{r}\,r^{\ell-1}\delta(r\eta)\,
            |f(k_{1}(\exp(r\eta))k_{2})|^{2}\notag\\[5pt]
    &=\,
        \int_{\mathfrak{a}^{+}}\diff{y^{+}}\,
            \delta(y^{+})\,|f(k_{1}(\exp{y^{+}})k_{2}|^{2}
    \label{last corollary estimate}
    \end{align}

Since $\diff{k_{1}}$ and $\diff{k_{2}}$ are normalized measures on $K$, we obtain finally
    \begin{align*}
        \|Sf\|_{L^{2}(G)}^{2}\,
    &=\,\const\,
        \int_{\mathfrak{a}^{+}}\diff{x^{+}}\,
        \delta(x^{+})\,
        |Sf(\exp{x^{+}})|^{2}\\[5pt]
    &\lesssim\,  
        \int_{\mathfrak{a}^{+}}\diff{x^{+}}\,
        \delta(x^{+})\,
        |\kappa_{1}(\exp{x^{+}})|^{2}\,
            \int_{K}\diff{k_{1}}\,\int_{K}\diff{k_{2}}\,
            |I_{k_{1},k_{2}}(|x^{+}|)|^{2}\\[5pt]
    &\lesssim\,
        \int_{K}\diff{k_{1}}\,\int_{K}\diff{k_{2}}\,
        \int_{0}^{\infty}\,\diff{R}\,R^{\ell-1}\,
        |I_{k_{1},k_{2}}(R)|^{2}\\[5pt]
    &\lesssim\,
        \int_{K}\diff{k_{1}}\,\int_{K}\diff{k_{2}}\,
        \int_{\mathfrak{a}^{+}}\diff{y^{+}}\,
            \delta(y^{+})\,|f(k_{1}(\exp{y^{+}})k_{2}|^{2}\\[5pt]
    &=\,\|f\|_{L^{2}(G)}^{2}.
    \end{align*}
Here, we have used successively the Cartan decomposition, the assumption of $\kappa_1$, the Jensen inequality, and estimate \eqref{last corollary estimate}.
\end{proof}

%%%%%%%%%%%%%%%%%%%%%%%%%%%%%%%%%%%%%%%%%%%%%%%%%%%%%%%%%%%%%%%%%%%%%%%%%%%%%%%%%%
%%%%%%%%%%                      ACKNOWLEDGEMENT                       %%%%%%%%%%%%
%%%%%%%%%%%%%%%%%%%%%%%%%%%%%%%%%%%%%%%%%%%%%%%%%%%%%%%%%%%%%%%%%%%%%%%%%%%%%%%%%%
\noindent\textbf{Acknowledgement.}
The authors thank helpful discussions with Jean-Philippe Anker. The authors are supported by the FWO Odysseus 1 grant G.0H94.18N: Analysis and Partial Differential Equations and by the Methusalem programme of the Ghent University Special Research Fund (BOF) (Grant number 01M01021). The first and second authors are supported by FWO Senior Research Grant G011522N. The second author is also supported by EPSRC grant EP/R003025/2.

%%%%%%%%%%% Bibliography
\printbibliography
%%%%%%%%%%% Bibliography

%%%%%%%%%%% Affiliation
\vspace{10pt}

\address{
    \noindent\textsc{Vishvesh Kumar:}
    \href{mailto:kumar.vishvesh@ugent.be}
    {vishvesh.kumar@ugent.be}\\
    \textsc{
    Department of Mathematics: 
    Analysis, Logic and Discrete Mathematics\\
    Ghent University, Belgium}
}\vspace{10pt}

\address{
    \noindent\textsc{Michael Ruzhansky:}
    \href{mailto:michael.ruzhansky@ugent.be}
    {michael.ruzhansky@ugent.be}\\
    \textsc{
    Department of Mathematics: 
    Analysis, Logic and Discrete Mathematics\\
    Ghent University, Belgium\\
    and\\
    School of Mathematical Sciences\\
    Queen Mary University of London, United Kingdom}
}\vspace{10pt}

\address{
    \noindent\textsc{Hong-Wei Zhang:}
    \href{mailto:hongwei.zhang@ugent.be}
    {hongwei.zhang@ugent.be}\\
    \textsc{
    Department of Mathematics: 
    Analysis, Logic and Discrete Mathematics\\
    Ghent University, Belgium}
}

%%%%%%%%%%%%%%%%%%%%%%%%%%%%%%%%%%%%%%%%%%%%%%%%%%%%%%%%%%%%
%%%%%%%%%%%%%%%%%%%%%%%%%%%%%%%%%%%%%%%%%%%%%%%%%%%%%%%%%%%%
%%%%%%%%%%%%%%%%%%%%%%%%%%%%%%%%%%%%%%%%%%%%%%%%%%%%%%%%%%%%
\end{document}